\documentclass[a4paper,11pt]{amsart}%{article}
\usepackage[margin=1.9cm]{geometry}
\usepackage{caption}
\usepackage[T1]{fontenc}
\usepackage{lmodern, amsfonts,amsmath,amstext,amsbsy,amssymb,
amsopn,amsthm,upref,eucal,mathptmx,mathtools,url}

\usepackage{mathrsfs}

\RequirePackage{xcolor} % [dvipsnames]
\definecolor{halfgray}
{gray}{0.55}%chapter numbers will be semi
\definecolor{webgreen}
{rgb}{0,0.4,0}
\definecolor{webbrown}
{rgb}{.8,0.1,0.1}
\definecolor{red}
{rgb}{1,0,0}
\usepackage{tikz}

\usepackage[hang]{subfigure}

\newcommand{\SL}{%
\operatorname{SL}
}

\newtheorem{theorem}{Theorem}[section]

\newtheorem{corollary}[theorem]{Corollary}

\definecolor{halfgray}
{gray}{0.55}%chapter numbers will be semi
\definecolor{webgreen}
{rgb}{0,0.4,0}
\definecolor{webbrown}
{rgb}{.8,0.1,0.1}
\definecolor{red}
{rgb}{1,0,0}
\usepackage{tikz}

\newcommand \R {{ \mathbb R}}

\date{\today}

\begin{document}

\begin{center}

\vspace{-9.9mm}

  {\bf \large {\sf XXI Congresso dell'Unione Matematica Italiana,
 Pavia, September 2019 }}
	\vspace{1.3mm}

  {\bf \large Slow chaos in surface flows}

	\vspace{1.9mm}
  {\bf \large Corinna Ulcigrai}
\end{center}

 \begin{abstract}
Flows on surfaces describe many systems of physical origin and are one of the most fundamental examples of dynamical systems, studied since Poincaré. In the last decade, there have been a lot of advances in our understanding of the chaotic properties of smooth area-preserving flows (a class which include locally Hamiltonian flows), thanks to the connection to Teichmueller dynamics and, very recentlty, to the influence of the work of Marina Ratner in homogeneous dynamics. We motivate and survey some of the recent breakthroughs on their mixing and spectral properties and the mechanisms, such as shearing, on which they are based, which exploit analytic, arithmetic and geometric techniques.
\end{abstract}

\section{Slowly chaotic dynamical systems}

\subsection*{Deterministic chaos and the butterfly effect.} 
\emph{Dynamical systems} provide mathematical models of systems which evolve in time. Many systems phenomena in our world, from the evolution of the weather to the motion of an electron in a metal, can be described by a dynamical system. While in a model one can include a \emph{random} component, or external \emph{noise}, we will restrict ourselves to fully \emph{deterministic} systems, whose evolution is completely described by a system of pre-determined rules or equations.  We will furthermore consider \emph{continuous time} dynamical systems, namely systems for which  the \emph{time} variable is a real parameter $t\in \mathbb{R}$,  described  by  a \emph{flow} on a space $X$, namely a $1$-parameter group\footnote{Assuming that $\varphi_{\mathbb{R}}=(\varphi_t)_{t\in\mathbb{R}}$ is a \emph{group} of diffeomorphisms under composition is equivalent to requiring  that  $\varphi_{t+s}(x)=\varphi_t(\varphi_s(x))$ for every $x\in X$ (or almost every in the measure-preserving set up introduced below) and every $t,s\in\mathbb{R}$. The typical example of a flow is given by solutions to differential equations. The precise definition of the type of  flows on which we will focus here, namely area-preserving flows on a surface, will be given  below.}  $\varphi_{\mathbb{R}}=(\varphi_t)_{t\in \mathbb{R}}$ of maps $\varphi_\mathbb{R}:X\to X$ (diffeomorphisms if $X$ is a smooth manifold).

Deterministic dynamical systems  often display \emph{chaotic features} (see Section \ref{sec:chaotic} for examples), which make  their behaviour as time grows hard to predict. This is a phenomenon known as \emph{deterministic chaos}. One of the best known features of chaotic behaviour is \emph{sensitive dependence on initial conditions} (SDIC for short), a property which was popularized as the \emph{butterfly effect}. In a system which displays SDIC, a small variation of the initial condition can lead to a macroscopically very different evolution after a long time. In particular, %if $x,y\in X$ are two points on a (metric) space $X$ and 
% for example,
 given a  flow $\varphi_{\mathbb{R}}: X\to X$ on a metric space $(X,d)$ with SFIC and  a point $x\in X$ (the \emph{initial condition}) one can find arbitrarily close initial conditions $y\in X$ such that the (forward) \emph{trajectories} of $x$ and $y$, namely  the orbits $(\varphi_t(x))_{t\geq 0}$ and $(\varphi_t(y))_{t\geq 0}$ drift apart\footnote{More precisely, if  $(X,d)$ is a metric space and   $\varphi_{\mathbb{R}}: X\to X$ a continuous flow,  $\varphi_{\mathbb{R}}$ has SDIC with sensitivity constant $\nu>0$ iff for every $x\in X$ and $\epsilon>0$ there exists $y$ such that $d(x,y)<\epsilon$ and there exists $t_0=t_0(\nu, y)$ such that $d(\varphi_{t_0}(x), \varphi_{t_0}(y))\geq \nu$.}. 

\subsection*{Fast chaos versus slow chaos.}
Dynamical systems can roughly be divided in three categories (\emph{hyperbolic}, \emph{elliptic} and \emph{parabolic}) according to the speed of divergence (if any) of close orbits.
%for example, on a  flow $\varphi: X\to X$ on a metric space $(X,d)$, given a point $x\in X$ (the \emph{initial condition}) one can find arbitrarily close initial conditions $y\in X$ such that the (forward) \emph{trajectories} of $x$ and $y$, namely  the orbits $(\varphi_t(x))_{t\geq 0}$ and $(\varphi_t(y))_{t\geq 0}$ drift further apart. 
In a \emph{hyperbolic} flow, the orbits $(\varphi_t(x))_{t\geq 0}$, $(\varphi_t(y))_{t\geq 0}$  of \emph{most}\footnote{\emph{Most}  means here \emph{almost every} with respect to a natural invariant measure for the flow. The formal definition of a (\emph{smooth}) hyperbolic flow (on the tangent space to a differentialble manifold) involve the existence of stable and unstable manifolds in the tangent space; the exceptional points $y$ whose trajectories do not diverge from (and actually converge to)   the trajectory of a given $x$ form the so-called stable manifold, which has positive codimension in the ambient space.} pairs $x,y$ of initial conditions diverge  \emph{exponentially} in time (i.e.~the distance $d(\varphi_t(x),\varphi_t(y))$, for small values\footnote{If the space $X$ if compact, the distance between two points has clearly an upped bound, so all statements about divergence must be interpreted \emph{infinitesimally}, i.e.~make sense for  of $t$ small, as asymptotic statements as the distance between $y$ and $x$ goes to zero, so that longer times can be considered.}  of time, is described by an exponential function of time. 
In a  \emph{parabolic} dynamical system, there is also \emph{divergence} of (most) nearby orbits, but this divergence happens at  \emph{subexponential} (usually \emph{polynomial}) speed. Finally while the flow is called \emph{elliptic} if there is no divergence  (or perhaps it is slower than polynomial).  Thus, both hyperbolic system and parabolic systems display SDIC, but the \emph{butterfly effect} happens
 at different speeds (respectivelly \emph{exponentially} or (sub)\emph{polynomially}).   We colloquially call these systems respectively \emph{fast chaotic} (when the butterfly effect is \emph{fast}, i.e.~exponential) and \emph{slowly chaotic} (when the butterfly effect is \emph{slow}, i.e.~polynomial or slower than polynomial).

While there is a classical and well-developed theory of hyperbolic systems (starting with the study of \emph{uniformely hyperbolic} dynamical systems, which was already developed in the $1970$s by mathematicians such as D.~Anosov and Y.~Sinai, Abel Prize in 2014, among others) and also a systematic study of  elliptic ones (starting with the theory of circle diffeomorphisms, whose study is intertwined with Hamiltonian dynamics and KAM theory), there is no general theory which describes the dynamics of parabolic flows and only classical and isolated examples are well-understood. 

\subsection*{Examples of parabolic systems.}\label{sec:examples}
Slowly chaotic (or \emph{parabolic}) systems include many dynamical systems of interest in physics, such as the \emph{Novikov model} of electrons in a metal (which will be discussed below), or the \emph{Ehrenfest model} (also called \emph{windtree model}) proposed by Paul and Tatjiana Ehrenfest in 1912 to undestand thermodynamics laws. 

Among examples arising in mathematics,  
perhaps the most studied (and better understoood) example of a parabolic flow is given in the context of hyperbolic geometry by the \emph{horocycle flow} on (the unit tangent bundle of) a compact negatively curved surface\footnote{In the context of \emph{homogeneous dynamics} (actions given by group multiplication on quotients of Lie groups), parabolic flows coincide with \emph{unipotent flows}. Horocycle flows  can be seen as the simplest example of unipotent flows on semi-simple Lie groups, given by the left action of upper triangular unipotent matrices in $\SL(2,\mathbb{R})$ on compact (or finite volume) quotients  $\SL(2,\mathbb{R})/\Gamma$.}:  
 while the \emph{geodesic} flow (whose trajectories move along geodesics for the hyperbolic metric) is a classical example of fast chaos and hyperbolic dynamics, when moving along \emph{horocycles} (which, in the upper half plane $\mathbb{H}$ are circles tangent to the real axis) one can show that divergence of nearbly trajectories is only a quadratic fuction of time, thus giving slowly chaotic dynamics.  

Another fundamental class of homogeneous  flows is given by \emph{nilflows}, or flows on (compact) quotients of \emph{nilpotent} Lie groups (\emph{nilmanifolds});  The \emph{prototype} example in this class are \emph{Heisenberg nilflows}, given by the action (by left multiplication) of a $1$ parameter subgroup of transformations of a compact quotient of the Heisenberg group\footnote{Let us recall that the Heisenberg group can be seen as the group $H$ of $3\times 3$ upper triangular matrices; a compact Heisenberg nilmanifold is obtained taking the quotient $H/\Gamma$ where $\Gamma<H$ is a discrete subgroup (for example the subgroup of matrices of the same form, but with \emph{integer} entries). A nilflow on a Heisenberg nilmanifold is then given by the action of a $1$-parameter subgroup $(h_t)_{t\in\mathbb{R}}\subset H$ by left multiplication $(g,t)\mapsto h_t g$.}.

In this survey we will focus on an another fundamental class of parabolic flows, in the context of \emph{area-preserving flows} on (higher genus) surfaces, %another important class of parabolic flows  is  %(by many authors, including us) considered 
given by \emph{locally Hamiltonian} flows, which are \emph{smooth} flows which preserve a \emph{smooth} area-form. The definition is given later on in \S~\ref{sec:locHam}.

Finally, starting from the classical examples of parabolic flows mentioned above, one can build new parabolic flows by considering \emph{perturbations}: the simplest perturbations are   \emph{time-changes} (or time-\emph{reparametrizations}) of a given flow, i.e. flows that move points along the \emph{same orbits}, but with different \emph{speed}.  More precisely, a flow $\tilde{\varphi}_\R$  is called  a (smooth)   
\emph{time-change}  of a flow  $\varphi_\R$ on $X$  if for all  
$x \in X$ and $t \in \mathbb{R}$ we have  $\varphi_t(x) = \varphi_{\tau(x,t)}(x) $ for some measurable
(smooth) function $\tau: X \times \mathbb{R} \rightarrow \mathbb{R}$. Notice that it follows from this definition that the time-change  $\tilde{\varphi}_\R$  has exactly the same trajectories than $\varphi_\R$ (but the motion along the trajectory has different speed).   Some time-changes, known as (smoothly) \emph{trivial}, give rise to flows that are (smoothly) \emph{conjugated} (i.e.~isomorphic as dynamical systems) to the original one and therefore have \emph{the same} chaotic properties. A feature of parabolic dynamical systems, though, is that 
 among smooth time-changes, {\it smoothly} trivial time-changes are \emph{rare}, i.e.~they often form a finite or countable codimension subspace\footnote{This is the consequence of the existence of \emph{distributional obstructions} (invariant distributions) to solve the cohomological equation. The first complete study of this phenomenon is perhaps Katok's work~\cite{Kat:CC1, Kat:CC2} on linear skew-shifts of the $2$-torus, which are closely related to Heisenberg nilflows. Let us also remark that \emph{finitely many} invariant distributions for horocycle flows in the finite area, non-compact case were first constructed by P.~Sarnak \cite{Sa:horo} by methods based on Eisenstein series. The structure of the space of obstructions was described in the case of translation flows (and locally Hamiltonian flows on surfaces) in \cite{Fo:sol}, for nilflows in~\cite{FF3} and for  horocycle flows in~\cite{FF1}.}. Therefore, the study of non trivial time-changes  allows to systematically produce new classes of parabolic flows.

\section{Chaotic properties}\label{sec:chaotic}
 A natural and fundamental question in parabolic dynamics (and dynamics in general) is which chaotic properties -in particular which fine ergodic and spectral properties- are \emph{generic} among  classes of (smooth) slowly chaotic flows. Let us give some examples in this section of which \emph{type} of \emph{chaotic properties} one might be interested in. We will comment on possible notions of \emph{generic} in \S~\ref{sec:generic}. 
 
Different type of chaotic properties are the focus of different branches of dynamics:  properties of topological nature (such as existence of dense trajectories) are studied in \emph{topological dynamics}, \emph{measure-theoretical} features (such as \emph{equidistribution} of a trajectory with respect to an equilibrum measure) in \emph{ergodic theory}, and properties of \emph{spectral nature} in \emph{spectral theory of dynamical systems}. 

\subsection*{{Topological dynamics}} focus on the most basic questions about the behaviour of trajectories, such as existence and abundance of \emph{periodic} trajectories (i.e.~trajectories of a point $x$ such that there exists a \emph{period} $t_0$ for which $\varphi_{t+t_0}(x)=\varphi_t(x)$ for all $t\in \R$), or existence and abundance of trajectories which are \emph{dense} in $X$. A (continuous) flow  $\varphi_\R:X\to X$ (where $X$ is a topological- or metric- space) is called \emph{minimal} if every orbit is \emph{dense}. In presence of fixed points (which is the case for surface flows in higher genus, which always have singularities for $g\geq 2$), the definition of minimality is slightly different: we only require all \emph{regular} orbits, namely all orbits $\varphi_\R(x)$ which are neither a fixed point nor a saddle separatrix,   to be  
dense.

\subsection*{{Ergodic theory}} studies flows which preserve a \emph{measure}: we assume hence that $(X,\mu)$ is a measure space and that $\varphi_\R:X\to X$ preserves the measure $\mu$, namely for any $A$ measurable set, $\mu(A)=\mu(\varphi_t(A)$ for all $t\in \mathbb{R}$. If a trajectory is \emph{dense}, one can further ask whether it is \emph{equidistributed} with respect to the invariant measure $\mu$, namely if the time spent in a measurable set $A$ is proportional (asymptotically) to its measure $\mu(A)$, or, in formulas, whether
\begin{equation}\label{eq:definition}
\lim_{T\to \infty}\frac{1}{T}\int_{0}^T \chi_A \left(\varphi_t(x)\right)\textrm{d}t = \mu(A).
\end{equation}
Systems for which this is true for \emph{almost every} initial condition $x$ (with respect to $\mu$) are \emph{ergodic}\footnote{Ergodicity is often defined in terms of \emph{metric indecomposability}: a measure preserving flow  $\varphi_\R:X\to X$  on  $(X,\mu)$ is \emph{ergodic} if any measureable $A$ which is \emph{invariant} under $\varphi_\R$, i.e. such that $\varphi_t(A)= A$ for all $t\in\R$, is measure-theoretically trivial, i.e. either $\mu(A)=0$ or $\mu (X\setminus A)=0$. The equivalence of this definition with the equidistribution property of almost every trajectory when $\mu(X)$ is finite is the content of the celebrated \emph{Birkhoff ergodic theorem}.}. A stronger conclusion, namely that equidistribution holds for \emph{every} point $x\in X$ with an infinite trajectoriy, holds for smooth flows which are \emph{uniquely ergodic}\footnote{Unique ergodicity is usually defined for a topological dynamical system, namely for \emph{continuous} maps or flows on a topological space $X$. We say that the system is \emph{uniquely ergodic} if there exists an \emph{unique} probability measure. One can show that in this case, for observables $f:X\to \mathbb{R}$ which are continuous, one gets a stronger conclusion that the Birkohff ergodic theorem, since one can show that  \eqref{eq:definition}, with $\chi_A$ replaced by $f$, holds for \emph{every} (and not only almost every) initial condition $x\in X$. This definition can be applied to the area-preserving flows which we consider later modulo some technicalities (in particular one has to exclude singularities) and guarantees that all trajectories of \emph{non-singular} points $x$ whose forward trajectory $(\varphi_t(x))_{t\geq 0}$ is not a separatrix are equidistributed in the sense that  \eqref{eq:definition} holds for all measurable sets $A\subset S$.}. 

A stronger property, \emph{mixing}, guarantees equidistributions not only of individual orbits, but of sets pushed under the flow $\varphi_\R$: in a mixing system, every measurable set $A\subset X$ \emph{equidistributes} (with respect to $\mu$) under the flow, i.e.
\begin{equation}\label{eq:mixingdef}
\lim_{t\to \infty}\mu (\varphi_t(A)\cap B) = \mu(A)\mu(B)
\end{equation}
for every measurable set $B$. This property is equivalent  to \emph{decay of correlations}, i.e. for every two smooth observables $f, g: X\to \R$, 
\begin{equation}\label{def:decaycorr}
\lim_{t\to \infty }\int_X \left( f\circ \varphi_t \right) g\, \mathrm{d}\mu - \int_X f \mathrm{d}\mu \int_X g\, \mathrm{d}\mu = 0. 
\end{equation}
This property is also known as \emph{strong mixing}, to distinguish it from another (weaker) property known as \emph{weak mixing} (where the convergence in \eqref{eq:mixingdef} is only required to happen along a subset of $t\in \R$ of \emph{density one}).  Mixing and weak mixing can also be interpreted as \emph{spectral properties},   (see footnote \ref{equivalentspectral}). Other type of mixing properties in addition to \emph{weak mixing} and \emph{strong mixing} include \emph{mild mixing} and \emph{mixing of all orders}. The latter generalizes mixing (which is defined using two sets $A,B$) to more sets: a  measure preserving flow $\varphi_\R$ on $(X,\mu)$  is  {mixing of order $N$}  if, for any $N$-tuple $A_0$,$\dots$, $A_{N-1}$ of measurables sets,
\begin{equation}\label{def:multiplemixing}
\mu \left(A_0 \cap \varphi_{{t_1}}({A_1})\cap \varphi_{{t_1+t_2}}({A_2}) \cap \cdots \cap \varphi_{{t_1+\dots +t_{N-1}}}({A_{N-1}}) \right)
 \xrightarrow{t_1,t_2\dots, t_{N-1} \to \infty}   \mu( { A_0}) \cdots \mu ({A_{N-1}} )
\end{equation}
and it is \emph{mixing of all orders} if  it is $N$-mixing for any $N\geq 2$. Equivalently, as for the definition of mixing, this can be reinterpreted as a statement about decay of \emph{multi}-correlations. It is a famous open conjecture, known as Rohlin's conjecture and still open, whether mixing implies mixing of all orders.

\subsection*{Spectral theory of dynamical systems} study the nature of the spectrum (and spectral measures) associated to the \emph{Koopman operator} a  (family of) operator(s) on $L^2(X,\mathscr{A},\mu)$  associated to measure preserving flow $\varphi_\R$ 
%namely they describe the nature of the spectrum of the \emph{Koopman} operator associated to a measure preserving dynamical system. This is the unitary opertor $U_t: L^2(X,\mu)\to L^2(X,\mu)$ on square-summable observables $f\in  L^2(X,\mu)$ obtained by pre-compositing with the flow, i.e. $f\mapsto U_t(f):= f\circ \varphi_t$.  
%(or transformation $T$) on the measure space $(X,\mathscr{A},\mu)$, namely the (family of) unitary operator(s) ($U^{\phi}_t$ or $U_T$) on $L^2(X,\mathscr{A},\mu)$ which acts on a square-integrable function $f$ by composing with %the dynamics, i.e. maps $f$ to 
(which acts by pre-composition $f\mapsto f\circ \varphi_t$ with the dynamics).  
One of the fundamental questions in spectral theory (see for example the surveys \cite{KT} or \cite{L} on spectral theory of dynamical systems)  is what is the \emph{nature of the spectrum} of the Koopman operator. To every $f\in L^2(X,\mu)$ one can associate  a  \emph{spectral} measure  denoted by  $\sigma_f$, i.e.\  the unique finite Borel measure on $\R$ whose Fourier coefficients are described by correlations, i.e.~such that
\begin{equation}\label{def:sp}
\int_X  f\circ \varphi_{t} \overline{f} \,d\mu =\int_\R e^{its}\,d\sigma_g(s)\quad\text{for every}\quad t\in \R.
\end{equation}
Spectral measures are useful to describe components of the the unitary representation given by the Koopman operator\footnote{Let us denote by $\R(g)\subset  L^2(X,\mu)$ the \emph{cyclic subspace} generated by $g$ which is given by $\R(g):=\overline{\langle \varphi_t(g):t\in\R \rangle} \subset  L^2(X,\mu)$. By the spectral theorem %(see e.g.~\cite{CFS:erg})
 the Koopman operator, restricted to $\R(g)$, 
is unitarily isomorphic to the $\R$-representation $(V_t)_{t\in\R}$ on $L^2(\R,\sigma_g)$ given by
$V_t(h)(s)=e^{its}h(s)$.}.  We say that the \emph{spectrum} of $\varphi_\R$ is (\emph{absolutely}) continous, or respectively (\emph{purely}) \emph{singular} iff for every $f\in L^2(X,\mu)$ the spectral measure $\sigma_f$ is (absolutely) continous, or respectively singular with respect to the Lebesgue measure on $\R$.
%(this \color{orange} is equivalent to the Koopman operator operator being \color{blue} equivalently implies that the \emph{spectral type} of the Koopman operator is \color{black} purely singular, see \cite{CFS:erg}).

Ergodicity, weak mixing and mixing can be expressed in terms of the spectrum of the Koopman operator\footnote{Ergodicity is equivalent to asking that the only eigenfunctions for the Koopman operator with eigenvalue $1$ are constant functions,\label{equivalentspectral} while weak mixing holds iff the Koopman operator has no non-constant eigenfunctions (i.e.~no eigenvalues different than $1$). Hence weak mixing implies that the spectrum is continous; mixing on the other hand is characterized by decay of the self-correlations $\langle U_t f, f\rangle_{L^2(X,\mu)}$, see \eqref{def:decaycorr}.}. %
%Mixing is equivalent to \emph{continuity} of the spectrum.}.  
Spectral results thus provide finer and stronger dynamical information.   For example, since mixing (and weak mixing), when they hold, provide, as spectral implication, the information that the spectrum is continuous (see footnote \ref{equivalentspectral}), proving that the spectrum is \emph{absolutely continuous} is e.~g.~ a strenghening of mixing, while \emph{singularity} of the spectrum shows that the system studied is far (more formally, \emph{spectrally disjoint}, a stronger concept of disjointess than that introduced in \S~\ref{sec:disjoint}) from strongly, fast chaotic systems. %The survey \cite{KT} provides an nice introduction and overview of questions and conjectures in spectral theory of dynamical systems for the interested reader. 
% We refer the interested reader to the survey \cite{KT} on spectral theory of dynamical systems. 

\subsection{Generic Chaotic properties in slowly chaotic systems}\label{sec:generic_chaotic}

There is a large and quite extensive literature on topological, ergodic and spectral properties of some of the  \emph{classical} parabolic examples mentioned in Section~\ref{sec:examples}. For example, the fine ergodic and spectral properties of the \emph{horocycle flow} have been studied in great detail\footnote{For \emph{compact} hyperbolic surfaces, it is for example known that the horocycle flow is minimal \cite{He}, uniquely ergodic~\cite{Fu:uni}, mixing~\cite{Pa:hor} (in fact it has Lebesgue spectrum \cite{St} and it is mixing of all orders \cite{Ma2}) and precise bounds  on  both the rate of mixing \cite{Rat:mix} and equidistribution % (for ergodic integrals of smooth 
%functions, see 
\cite{FF1} are available. For the case of finite-volume, non-compact surfaces, see e.g.~the works  by  Dani \cite{dani} and  %the horocyle flow is not uniquely ergodic
%and the classification of invariant measures is due to Dani 
 %The asymptotic behaviour of averages along closed horocycles has been studied  by D.~Zagier~\cite{Za}, P.~Sarnak~\cite{Sa}, D.~Hejal \cite{Hj} and more recently in \cite{FlaFo} and by A.~Str{\"o}mbergsson \cite{St}. Horocycle flows on general geometrically finite surfaces have been studied by 
M.~Burger~\cite{Bur}. 
} and mostly were already well understood in the $1970s$ (see for example \cite{Fu:uni, Ma2, Pa:hor,  FU, Bur, FF1, Str} and more in general \cite{St} or \cite{AGH}, and the reference therein, for unipotent flows). 

It is well known that the typical\footnote{\emph{Typical} means here for a full measure set of frequencies of the underlying toral factor, defined as follows. Every nilmanifold $G/\Gamma$ is a fiber bundle over a torus. In fact, the group~$\overline{G}=[G,G]\backslash G$ is Abelian, connected and simply connected, hence isomorphic to~$\R^n$ and~$\overline{ \Gamma}=[\Gamma,\Gamma] \backslash \Gamma$ is a lattice in~$\overline{G}$. Thus we have a natural projection $p: G/\Gamma \to \overline{G}/ \overline{\Gamma}$ over a torus of dimension~$n$. The nilflow hence project to a linear flow on an $n$-dimensional torus. We say that a property holds for a   \emph{typical}  nilflows if it holds for a \emph{typical} linear flow on the underlying toral factor, namely for a full measure set of rotation numbers (also called \emph{frequencies}).\label{toralfactor}} nilflow is minimal and uniquely ergodic \cite{AGH}, however, nilflows are never (weak) mixing, due to an intrinsic obstruction, namely the presence of a toral factor\footnote{As explained in the foonote \ref{toralfactor} above, \label{nilflownevermixing} any nilflow project on a linear flow on a toral factor.   It follows that a nilflow is never weakly mixing (and hence never mixing), since the linear toral flow has pure point spectrum and hence many non-trivial eigenfunctions (the toral characters), which can be pulled back to the nilmanifold to produce eigenfunctions for the nilflow (recall that weak mixing can be characterized in terms of lack of eigenfunctions , see \ref{sec:chaotic}). However, it is possible to prove by methods of representation theory that any nilflow is  {\it relatively mixing}, in the sense that the limit of correlations of functions with zero average along all fibers of the projection $p$ is equal to zero. Nilflows also 
have the property of \emph{relative} Lebesgue spectrum (namely, the spectrum restricted to observables in the 
orthogonal complement of the span of the pull-back  of the toral characters, is countable Lebesgue spectrum, see \cite{AGH}).}. 
Results on the speed of equidistribution of Heisenberg nilflows for smooth functions 
were proved by L.~Flaminio and Forni in \cite{FF2}. % who showed that  
%In fact, it is possible to prove 
%by the theory of unitary representations of the Heisenberg group (the Stone-Von Neumann theorem, see for example \cite{CG:rep}, \S 2.2)  that for all 
%for sufficiently smooth functions in $H^\perp$ the decay of correlations is polynomial (it
%is faster than any polynomial for infinitely differentiable functions in $H^\perp$).
%rates of equidistribution were
A series of recent works \cite{AFU, Rav:nil, AFRU} indicates that, even though classical nilflows are never  mixing (see footnote \ref{nilflownevermixing}), a typical \emph{time-change} (in a dense class of smooth time-changes) of a mimimal nilflow on any nilmanifold (different from a torus) is mixing.

In the rest of this survey  we will focus  on generic chaotic properties of smooth area preserving flows. While the understanding of minimality and ergodicity follows from results from the $1970s$ and $1980s$ respectively, the study of mixing properties has been the object of active research in the last decade or so, while the first breakthroughs on spectral properties are only very recent.

%; \cite{Ka:int, Keane, Ma:int, Ve:gau, AF:wea} or lecture notes such as \cite{Ma:rat, Yo, FoMa} and the references therein for translation flows;  \cite{Ul:wea, Ul2, Rav1, KKU, BK} among others for locally Hamiltonian flows.  

These results do not show an entirely coherent picture and make the identification and description of characteristic \emph{parabolic} features difficult. For example, while the horocycle flow (as well as all its smooth time-changes) are mixing (and actually mixing of all orders, see \cite{Ma1, Ma2}), as we recalled above (see in particular footnote \ref{nilflownevermixing}) nilflows are never mixing. This  difference in behavior can be attributed to the lack of parabolicity in certain directions (those that live in  the toral factor, see foonote \ref{toralfactor}), which would suggest calling  nilflows \emph{partially parabolic} systems. Nevertheless, this obstruction %  In this paper we prove that these obstructions
   can be broken by a perturbation (as we show in \cite{AFRU}, see also \cite{AFU} for the special case of Heisenberg nilflows), so  that in a dense set of smooth-time changes all flows which are not trivially conjugate to the nilflow itself are indeed mixing. Similarly, recent results seem to indicate that certain disjointess properties (see \S~\ref{sec:disjoint}), which do not hold for a classical example such as the horocycle flow, are nevertheless generic among time-changes. 
	
This shows that  to better understand slowly chaotic systems it is therefore crucial both to understand what are the finer chaotic properties of the known parabolic examples as well as,  at the same time, to study new classes of parabolic examples (such as those produced by time-changes and other parabolic perturbations\footnote{It should be pointed out that describing more general perturbations (beyond the class of time changes) which produce (new) \emph{parabolic} flows is quite delicate,  since by a perturbation  one \emph{typically} gets a hyperbolic flow. Examples of parabolic perturbations, which are not time changes, can be constructed for example by \emph{twisting} (see for example the work of Simonelli on \emph{twisted} horocycle flows, in~\cite{Si}).  New examples of parabolic perturbations for which one can study ergodic theoretical properties were recently constructed by Ravotti~in~\cite{Rav3}.}. 
    % mixing is a generic property among parabolic perturbations of nilflows.

\section{Chaotic properties of smooth surface flows.}\label{sec:surfaces}
We will focus now on the class of slowly chaotic systems given by smooth  area-preserving  (or locally Hamiltonian) flows on surfaces, surveying the recent advances in our understanding of their mixing, spectral and disjointness  properties as well as the mechanisms which explain them.

\subsection{Flows on surfaces.}\label{sec:flowssurfaces}
Flows on surfaces are one of the most basic and most fundamental examples of dynamical systems, whose study goes back to Poincar{\'e} \cite{P} at the end of the Ninetineth centrury,  and coincides with the birth of dynamical systems as a research field. 
Many models of systems of physical origin are described by surface flows:  Poincar{\'e} motivation to study surface flows was for example related to celestial mechanics and the two physical systems already mentioned before,  the Ehrenfest model in statistical mechanis  and the Novikov model in solid state physics, can be described by flows on surfaces (respectively linear flow on an translation surface and to a locally Hamiltonian flows). 

In addition to providing a  fundamental classes of \emph{parabolic} dynamical systems, smooth, area-preserving flows on surfaces, are fundamental in dynamics  because they are  among the lowest possible dimensional smooth dynamical system (on compact manifolds of lower dimension, the  other fundamental class of smooth dynamical systems are circle diffeomorphisms, whose rich theory is a cornerstone of dynamics). Despite having zero entropy, as shown in \cite{LSY}, they neverthelss display a rich variety of chaotic properties and, despite their basic nature, there are still many open questions on the mathematical characterization of chaos (in particular on dynamical, spectral and rigidity question) in various natural classes of surface flows. 

In this survey we will only be  concerned with flows which preserve a (probability) measure (see \S~\ref{sec:chaotic} for the definition), for example an area-form, since this is the natural setup for \emph{ergodic theory} (see Section \ref{sec:chaotic}). % (the branch of dynamics which investigates ergodic properties of measure-preserving dynamical systems).  %Let us recall that a flow  $(\varphi_t)=(\varphi_t)_{t\in \mathbb{R}}$ is a $1$-parameter groups of diffeomorphisms $\varphi_t: S \to S $. We will assume that $(\varphi_t)$ preserves a measure $\mu$ (namely $\mu(\varphi_t (A))=\mu(A)$ for any measurable set $A$ and $t\in \mathbb{R}$).

\subsection{Locally Hamiltonian flows}\label{sec:locHam}
Let $S$ be a compact, connected, orientable (smooth) surface and let $g$ denote its genus.  We will assume throughout that $g\geq 1$. Perhaps the most natural class of measure-preserving flows on $S$ are smooth flows preserving a smooth measure (with smooth absolutely continuous density).   Let  $\omega$ be a fixed smooth area form (locally given in coordinates $(x,y)$ by $f(x,y) dx\wedge dy$ where $f$ is a smooth function). Thus, equivalently, the pair $(S, \omega)$ is a two-dimensional symplectic manifold. We will consider a smooth flow $\varphi_\R=(\varphi_t)_{t\in\mathbb{R}}$ on $S$ which preserves a measure $\mu$ given  integrating a smooth density with respect to $\omega$. We will assume that the area is normalized so that  $\mu(S) =1$. 
It turns out that such smooth area preserving flows on $S$  are in one-to-one correspondence  with smooth closed real-valued differential $1$-forms as follows.  Given a smooth, closed, real-valued differential $1$-form $\eta$, let $X$ be the vector field determined by $\eta = i_X \omega$ where $i_X$ denotes the contraction operator, i.e. $i_X \omega =\omega( \eta, \cdot )$ and consider the flow $\varphi_\R$ on $S$ given by $X$. Since $\eta$ is closed, the transformations $\varphi_t$, $t \in \mathbb{R}$, are  area-preserving. Conversely, every smooth area-preserving flow can be obtained in this way.

 \begin{figure}[h!]
 \subfigure[A flow on a surface of $g=3$\label{g3}]{  \includegraphics[width=0.5\textwidth]{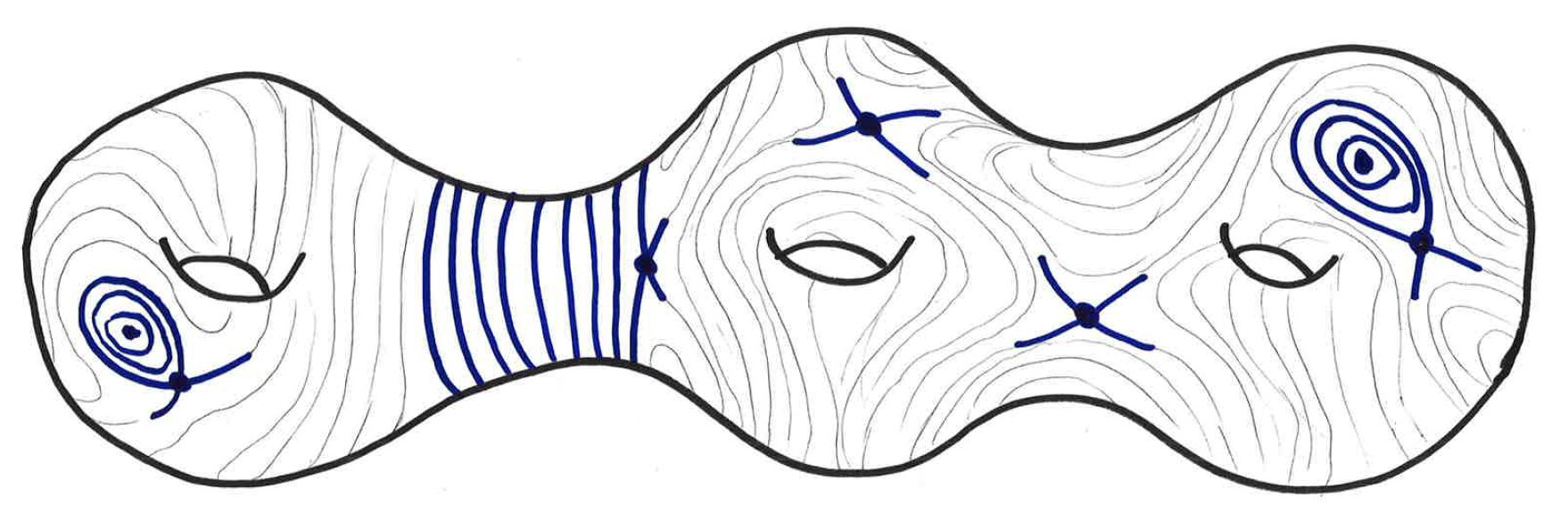}}\hspace{9mm}
 \subfigure[An Arnold flow ($g=1) $\label{Arnoldtorus}]{
  \includegraphics[width=0.27\textwidth]{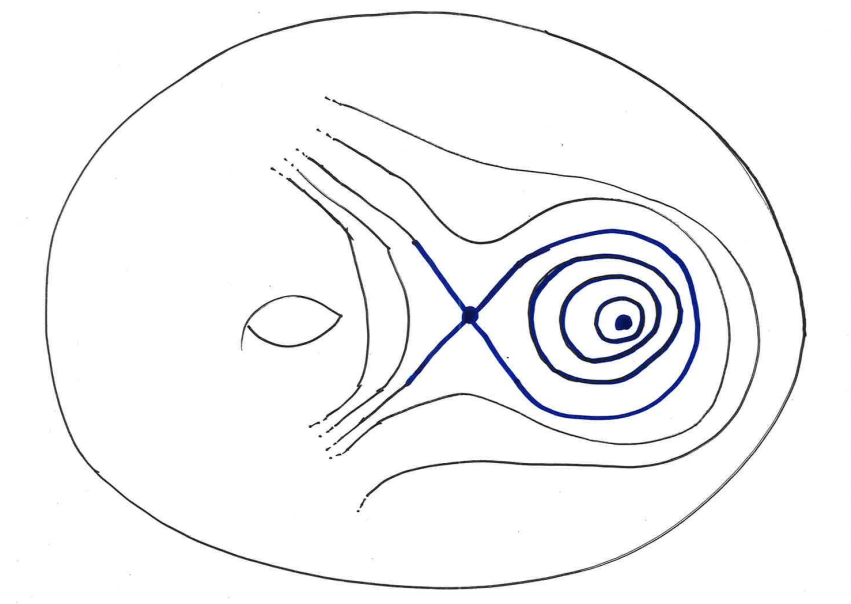} 	} 
%\subfigure[on $\mathbb{R}^2/\mathbb{Z}^2$ \label{Arnoldflowsquare}]{\includegraphics[width=0.33\textwidth]{LocHamSquare}   }	
 \caption{Trajectories of locally Hamiltonian flows on a surfaces. \label{locHamflows}}
\end{figure}

%\begin{figure}[h!]
%\caption{Decomposition in periodic components filled by closed orbits (two islands around centers and a cylinder in blue in the Figure) and minimal components (one of of genus one and one of genus  two in the example).\label{decomposition}}
%\end{figure}

The flow $\varphi_\R$  is known as the \emph{multi-valued Hamiltonian} flow associated to $\eta$. Indeed, the flow $\varphi_\R$ is \emph{locally Hamiltonian}, i.e.\ \emph{locally} one can find coordinates $(x,y)$ on $S$ in which $\varphi_\R$ is given by
%The form $\eta$ determines a flow  $\varphi_\R$ %given by a \emph{multi-valued Hamiltonian} as follows. 
 % a closed surface of genus $g\geq 2$ with a fixed area form and a closed differential $1$-form $\eta$ on it. 
%since $\eta $ is closed, one can locally write $\eta= \ud H$ for some real-valued function $H$ and 
 the solution to the  equations $$\begin{cases}\dot{x}&={\partial H}/{\partial y},\\ \dot{y}& =-{\partial H}/{\partial x}\end{cases}$$ for some smooth  real-valued Hamiltonian function $H$.  A \emph{global}  Hamiltonian $ H$ cannot be in general be defined (see \cite{NZ:flo}, Section 1.3.4), but one can think of  $\varphi_\R$ as globally given by a \emph{multi-valued} Hamiltonian function.

 \begin{figure}[h!]
  \subfigure[center \label{center}]{
  \includegraphics[width=0.23\textwidth]{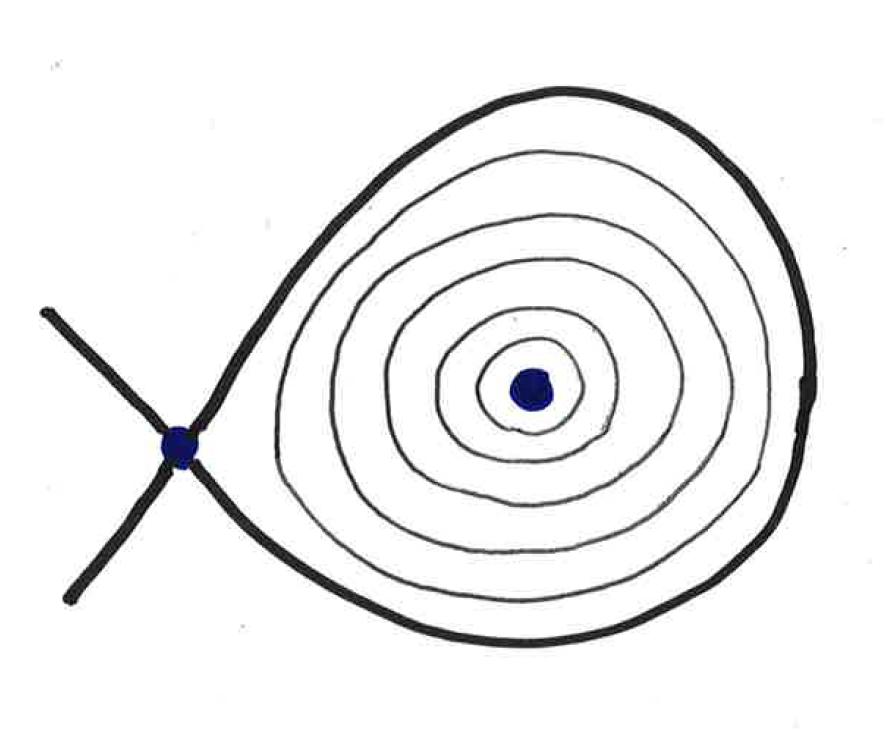} 	} \hspace{6mm} \subfigure[simple saddle \label{simplesaddle}]{ \includegraphics[width=0.18\textwidth]{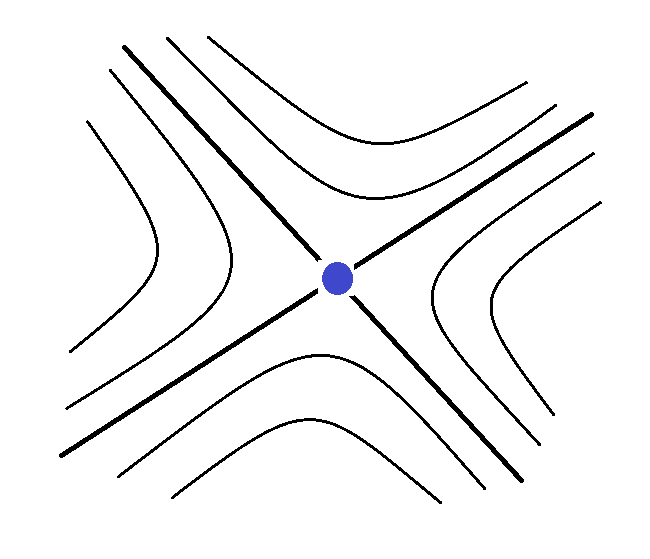}}\hspace{9mm}
\subfigure[\label{multisaddle}]{
\includegraphics[width=0.18\textwidth]{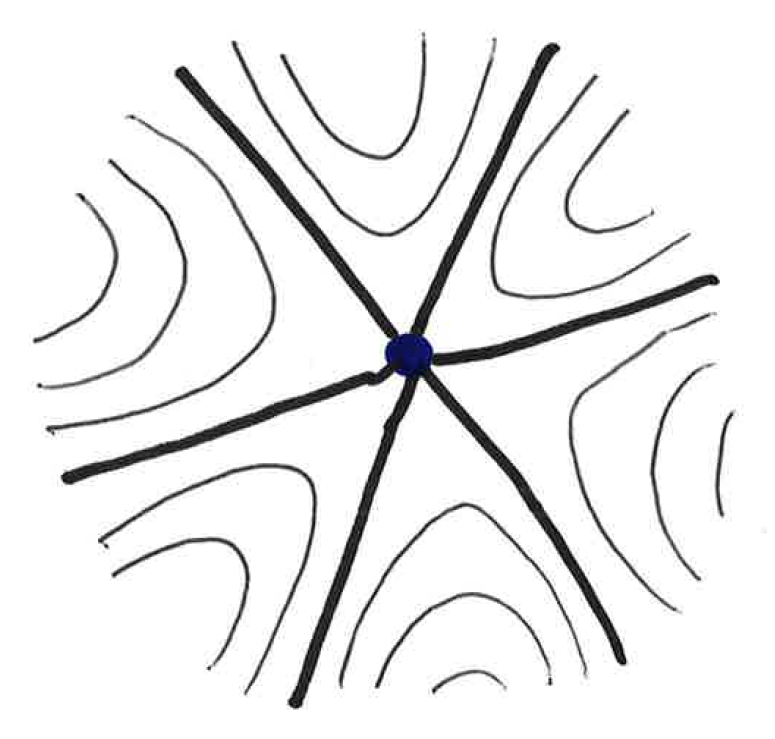}	\hspace{9mm}  }
{\caption{Type of singularities (i.e.~fixed points) of a locally Hamiltonian flow.\label{sing_types}}}
\end{figure}

%Exaples of trajectories 
Let us remark that locally Hamiltonian flows necessarily  have fixed points, or \emph{singularities}, if $g \geq 2$. Singularities, as shown in Figure~\ref{sing_types}, can be either centers (Fig.~\ref{center}), simple saddles (Fig.~\ref{simplesaddle}) or multi-saddles (i.e.~saddles with $2k$ pronges, $k\geq 2$, see Fig.~\ref{multisaddle} for $k=3$), . Examples of flow trajectories are shown in Figure~\ref{locHamflows}.  For $g=1$, i.e.~on a torus, if there is a singularity than there has to be another one. The simplest examples of locally Hamiltonian flows  with singularities on a torus, i.e.~flows with one center and one simple saddle (see Figure~\ref{Arnoldtorus}), were studied by V.~Arnold in \cite{Ar:top} and are nowadays often called \emph{Arnold flows}\footnote{More precisely, referring to the decomposition described in \S~\ref{sec:minimal}, we call \emph{Arnold flow} the restriction to a minimal component obtained by removing the center and the disk filled by periodic orbits around it (called \emph{island}), which, as Arnold shows in  \cite{Ar:top}, is always bounded by a saddle loop.\label{Arnoldflowfootnote}.}

%\begin{figure}[h!]
%{    \includegraphics[width=0.5\textwidth]{ArnoldSpecialFlow}}
%%\centering{\caption{
%Special representation of an {Arnol'd} flow. \label{Arnoldspecialflow}}}
%\end{figure}

% \begin{figure}[h!]
  %\end{figure}

A lot of interest in the study of multi-valued Hamiltonians and the associated flows - in particular, in their ergodic and mixing properties - was sparked 
by Novikov \cite{No} in connection with problems arising in solid-state
physics as well as in pseudo-periodic topology (see e.g.\ the survey \cite{Zo:how} by A.~Zorich). Indeed, Novikov \cite{No} and his school in the $1990s$  advocated the study of locally Hamiltonian flows as model to describe the motion of an electron in a metal under a magnetic field in the semi-classical approximation (the surface appears here as Fermi energy level surface). Novikov made some conjectures (known as \emph{Novikov problem}) on the asymptotic behaviour of trajectories of electrons. At the same time, Arnold \cite{Ar:top} made a conjecture on mixing for the flows  we call today \emph{Arnold flows} (see footnote \ref{Arnoldflowfootnote}). This conjecture has been the motivation for a lot of the work on the mixing properties of locally Hamiltonian flows, see the overview given in \S~\ref{sec:mixing}.

In order to survey the current knowledge of chaotic properties of locally Hamiltonian flows, it is useful to first point out their relation with another well studied class of  area-preserving flows on surfaces, namely linear flows (which, we stress for the reader, are not \emph{smooth} surface flows).

\subsection{Linear flows and time-changes of locally Hamiltonian flows.}\label{sec:linear} 
 The basic example of a \emph{linear flow} is the flow given on the torus $\mathbb{R}^2/\mathbb{Z}^2$ by solutions of 
 $$\left( x'(t) , y'(t)\right)=\left( \cos \theta,   \sin \theta\right),$$
  which move along at unit speed along  (the image in $\mathbb{R}^2/\mathbb{Z}^2$ of) Euclidean lines. Linear flows (also called \emph{translation flows}) can be defined more in general on \emph{translation surfaces}, namely surfaces which are locally Euclidean outside a finite number of conical singularities (with cone angles $2\pi k, k\in \mathbb{N}$, which produce saddles of the flows with $2k$ prongs, see Figure~\ref{multisaddle} for $k=3$). Notice that these flows preserve a Euclidean area (but are discontinuous flows, since singularities are reached in finite time).

It turns out that every \emph{minimal} locally Hamiltonian flow on $S$ (as well as the restriction of a locally Hamiltonian flow to one of its 
minimal components, see \S~\ref{sec:minimal}), in suitably chosen coordinates, are time-changes (or time-reparametrizations, we refer to \S~\ref{sec:examples} for the definition) of a linear flow on $S$ (or a subsurface of $S$ in the case of a minimal component). Thus, minimal locally Hamiltonian flows 
 have the same \emph{trajectories} (up to time-reparametrization) than  \emph{linear flows} on translation surfaces (see for example~\cite{Zo:how}). In particular, certain properties, such as minimality and ergodicity (as well as homological aspects such as the asymptotic behaviour in the Novikov problem), which depend only on trajectories as sets and not on their time-parametrization, can be deduced  for locally Hamiltonian flows by studying them in linear flows.   This was in part one of the original motivations (in addition to unfolding of rational billiards in the West) that sparked the interest of mathematicians such as A.~Zorich in the ergodic theory of \emph{linear flows}  (see below). 
 %and led to the blossoming research in this related area.  
 We stress though that,   while some chaotic properties like ergodicity  depend on the orbits of the flow, others (like mixing and spectral properties) crucially depend on the time-parametrization of the orbits and require ad-hoc techniques (see \S\ref{sec:shearing}and \S\ref{sec:spectral}). 

%The study of linear flows is intertwined with that of \emph{interval exchange transformations} (or in brief IETs). These are piecewise isometries of the unit interval, which were introduced by Oseledets in the $1960s$ as a basic example to study in ergodic theory, and also appear naturally as Poincar{\'e} (or first return) map of linear flows on a section. More formally, a IET $T:[0,1]\to [0,1] $ is a bijection which acts as a rigid translation on each interval $I_i$ (called \emph{continuity interval}) of a finite partition such that $I_1 \cup\dots \cup I_d = [0,1]$.  
%A $d$-IET, or a IET with $d$ continuity intervals, is completely determined by a permutation $\pi$ on $d$-symbols and the lenghts $\lambda_i, 1\leq i \leq d$ of the $d$ continuity intervals. 
%We will say that a property holds for \emph{almost every} IET if it holds for Lebesgue almost every choice of  the lenghts $Leb(I_i)$,  $ 1\leq i \leq d$ of the 
%continuity intervals %.  $(\lambda_1, \dots, \lambda_d)$ 
%(under a natural condition on the permutation  $\pi$ known as \emph{irreducibility}). %%% Possibly to shorten

\subsection{Linear flows and Teichmueller dynamics.}\label{sec:Teich} 
The study of linear flows on translation flows  and their ergodic properties has been a highly topical area of research for the past four decades (from the $1980$s), in connection with the study of billiards in (rational) polygons, interval exchange transformations (or for short IETs) and Teichmueller dynamics, a research area  which has benefited from the contribution of several Fields medallists (including Avila, Kontsevich, McMullen, Mirzakhani and Yoccoz).  %{\color{blue} Mention programs (in Bonn at the Hausdorff Insitute and  in 2012, Fields 2019, MSRI 2020)?}
%Let us briefly summarize the key achievements. 

In virtue of this flourishing activity, the ergodic and spectral properties of \emph{typical} (in the measure theoretical sense)  translation flows  are  by now well understood. Let us say that a property holds for a \emph{typical} linear flow if it holds in a.e.~direction on  \emph{almost every} translation surface with respect to a natural measure on translation surfaces known \emph{Masur-Veech measure}\footnote{Perhaps the simplest way to define the Masur-Veech measure, is to consider a presentation of a translation surface $S$ as a polygon with $2N$ sides with  pairs of parallel, congruent sides $(v_1,v_1'), \dots, (v_N, v_N')$ identifyed by glueings. Then  the vectors $(v_1, \dots, v_N)\in \mathbb{R}^N$ give local coordinates for an open set $U$ of translation surfaces around $S$ and the Masur-Veech measure is just Lebesgue measure on $U\subset \mathbb{R}^N$.}. One of the first results,  shown 
 already in the $1970s$ by M.~Keane \cite{Keane}, is that a typical linear flow is minimal. Moreover, it is  \emph{uniquely ergodic} (which implies that every infinite orbit is not only dense, but  also equidistributed, see \S~\ref{sec:chaotic}). Unique ergodicity of typical linear flows was known as \emph{Keane's conjecture} and 
  proved independently in the seminal works by Masur and Veech \cite{Ma:int, Ve:gau} through renormalization techniques which gave birth to the topical field of  Teichmueller dynamics. On the other hand  linear flows are never mixing, as proved by Katok \cite{Ka:int} already in  the $1980s$,  but they are neverthelss typically (in the above sense) \emph{weak-mixing} (refer to Section~\ref{sec:chaotic} for definitions).
a long-standing conjecture settled by Avila and Forni in  \cite{AF}. % or the \emph{Kontsevich-Zorich conjecture} \cite{Zo:how}

From the spectral theory perspective, for \emph{typical} translation flows, the nature of the spectrum (which turns out to be singular continuous) has been known since seminal work by Veech, see \cite{Ve:gau, Ve:IET}).  
%on \emph{weak-mixing} of almost every IET \cite{AF}
%see e.g. the proof that almost every IET (and almost every linear flow) is \emph{weak-mixing} by Avila and Forni  \cite{AF} %\cite{AF:wea} %%$%Shortableor the 
%proved 
%settled by Avila and Viana \cite{AV}. %\cite{AV:Lya} 
%of simplicity of the Lyapunov spectrum, that settled 
%of the deviations phenomenon described in the \emph{Kontsevich-Zorich conjecture} \cite{Zo:how}. %this is a phenomenon whose discovery was 
%originally motivated by the understanding of Novikov flows, see \cite{Zo:how, AV2}). 
Recently there have been also advances in the spectral theory of non generic (especially self-similar) translation flows and IETs,  see for example %the works by Bufetov and Solomyak
 \cite{BSU, BS1, BS2}.  %on spectral measures of IETs and their modulus of continuity. 
%While minimality and ergodicity for (classes of) locally Hamiltonian flows follow directly by Keane and Masur and Veech work (see below), mixing and spectral properties for locally Hamiltonian flows have to b

\smallskip
We now discuss locally Hamiltonian flows, explaining how minimality and unique ergodicity can be understood reducing the study of (minimal components of) locally Hamiltonian flows to linear flows, while the classification of mixing properties (described in \S~\ref{sec:mixing}) has been based on geometric mixing mechanisms specific to smooth slowly parabolic flows (such as \emph{shearing}, see \S~\ref{sec:shearing}) and the spectral theory is only now starting to be understood (refer to \S~\ref{sec:spectral}).  We first specify (in the next \S~\ref{sec:generic}) the topology and measure that we will use on the space of locally Hamiltonian flows.

\subsection{Genericity notions for locally Hamiltonian flows} \label{sec:generic}
Let us define two natural ways of definining a notion of \emph{generic} (or \emph{typical}) locally Hamiltonian flow, one topological and the other measure-theoretical.

\smallskip
One can define a \emph{topology} on locally Hamiltonian flows by considering perturbations of closed smooth $1$-forms by (small) closed  smooth $1$-forms\footnote{Let $\eta$, $\eta'$ be two smooth closed $1$-forms. We say
that $\eta'$ is an $\epsilon$-perturbation of $\eta$ if for any $x\in S$ there exists coordinates on a simply connected neighbourhood  $U$ of $x$, such that  $\eta \Vert _U=dH $ and $(\eta'-\eta)\Vert_U= dh$ where $\Vert h\Vert_{\infty}<\epsilon \Vert H\Vert_{\infty}$ (here $\Vert \cdot \Vert_{\infty}$ denotes the $\mathscr{C}^{\infty}$ norm).}. We say that a condition is \emph{generic} (in the sense of Baire) if it holds for
flows described by an open and dense set of forms with respect to this topology. For example, asking that the $1$-form $\eta$  is \emph{Morse}, i.e.\ it is locally the differential of a Morse function (which has non-degenerate zeros) is a generic condition.
%  Thus, all zeros of $\eta$ correspond to either centers or simple saddles.   This condition is generic (in the Baire cathegory sense) in the space of perturbations of closed smooth $1$-forms by closed smooth $1$-forms. 

\smallskip
 A \emph{measure-theoretical notion of  typical} is defined as  follows by using the \emph{Katok fundamental class} (introduced by Katok in \cite{Ka:inv}, see also \cite{NZ:flo}), i.e.\ the cohomology class of the 1-form $\eta$  which defines the flow.  Let $\Sigma$ be the set of fixed points of $\eta$  and let $k$ be the cardinality of $\Sigma$. Let $\gamma_1, \dots, \gamma_n$ be a base of the relative homology $H_1(S, \Sigma, \mathbb{R})$, where $n=2g+k-1$. The image of  $\eta$ by the period map $Per $ is $Per(\eta) = (\int_{\gamma_1} \eta, \dots, \int_{\gamma_n} \eta) \in \mathbb{R}^{n}$. The pull-back $Per_* Leb$ of the Lebesgue measure class by the period map gives the desired measure class on closed $1$-forms. When we use the expression \emph{typical} below, we mean full measure with respect to this measure class.

\subsection{Periodic and minimal components}\label{sec:minimal}
\begin{figure}[h!]
 \subfigure[an island of period orbits \label{island}]{
 \includegraphics[width=0.2\textwidth]{Island} 	} \hspace{6mm}
 \subfigure[a cylinder of periodic orbits \label{cylinder}]{
	\includegraphics[width=0.2\textwidth]{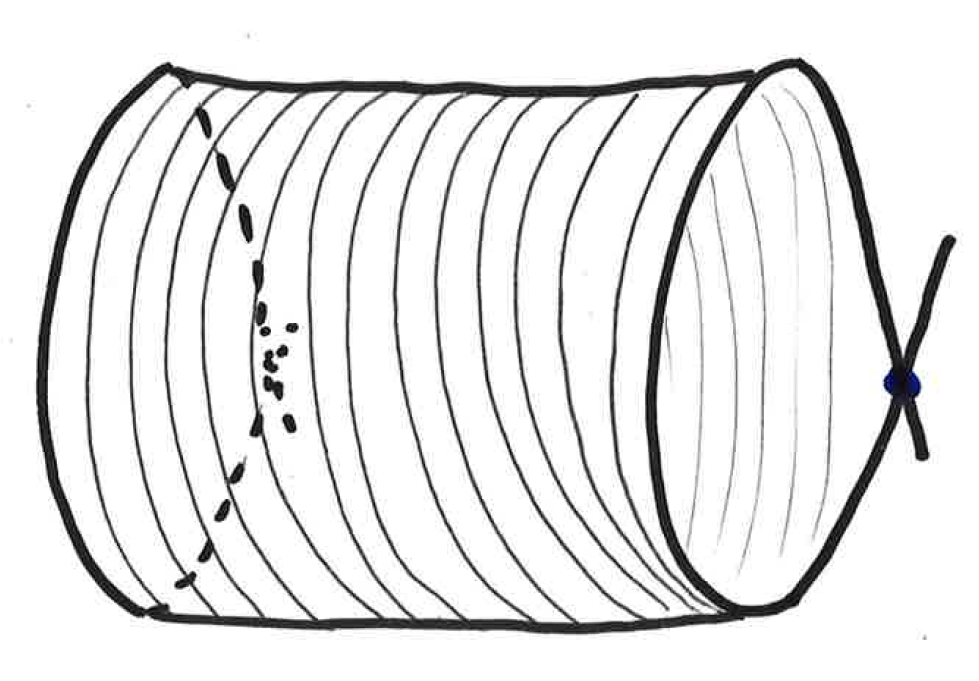}} \hspace{6mm}
\subfigure[a minimal component of $g=2$ \label{minimalcomp}]{\includegraphics[width=0.33\textwidth]{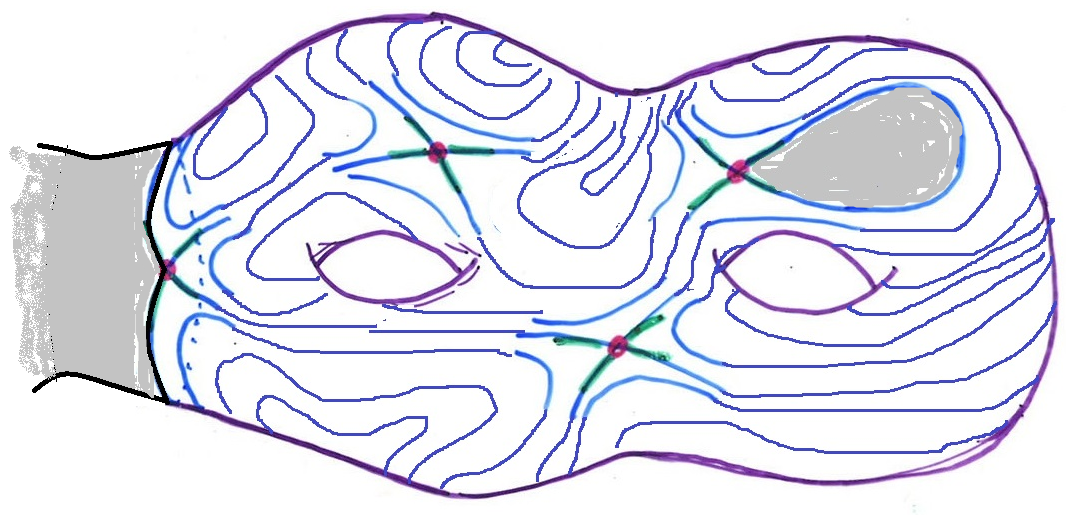}   }	
 \caption{Examples of periodic and minimal components% in the decomposion of a locally Hamiltonian flow
.\label{components}}
\end{figure}

%Recall When $g\geq 2$, the  (finite) set of fixed points (or singularities) of $\varphi_\R$ is always non-empty.
  A \emph{generic} locally Hamiltonian flow (in the sense of Baire category, with respect to the topology defined in the previous \S~\ref{sec:generic})  has only \emph{non-degenerate fixed points}, i.e.\ \emph{centers}  and \emph{simple saddles} (see (see Figures~\ref{center} and ~\ref{simplesaddle}), as opposed to degenerate \emph{multi-saddles} (as in Figure \ref{multisaddle}).  We call  \emph{saddle connection} a flow trajectory from a saddle to a saddle and a \emph{saddle loop} a saddle connection from a saddle to the same saddle (see Figure \ref{island}). It can be shown that each center is contained in a disk filled with closed (i.e.~periodic) trajectories and bounded by a saddle loop, called an \emph{island} of periodic orbits, see Fig.~\ref{island}. Hence, in presence of centers, the flow $\varphi_\R$ is never minimal (since orbits in the complement of the island avoid the island and hence cannot be dense).

	From the point of view of topological dynamics (as proved independently  by Maier \cite{Ma:tra}, Levitt \cite{Le:feu} and Zorich \cite{Zo:how}), each smooth area-preserving flow can be decomposed into subsurfaces (with boundary) on which the restriction of $\varphi_\R$ either foliates into closed (i.e.~periodic) orbits and up to $g$ subsurfaces (recall that $g$ is the genus of $S$) 
on which (the restriction of) $\varphi_\R$ 
 is \emph{minimal}, i.e.~every bi-infinite orbit is dense. The first ones are called \emph{periodic components} and are either  islands (as in Fig.~\ref{island}) or cylinders filled by periodic orbits and bounded by saddle loops, as in Fig.~\ref{cylinder}. The latter are known as \emph{minimal componens}  (see an example in Fig.~\ref{minimalcomp}) and by topological reasons there cannot be more than $g$ of them. The flows in Figure~\ref{locHamflows}, for example, can be decomposed, in the case of \ref{g3}, into two islands and one cylinder filled by closed orbits and two minimal components  (one of of genus one and one of genus two), while, in the case of the flow on the torus in Fig.~\ref{Arnoldtorus}, there is one island and one minimal component (the so-called Arnold flow). 
%( Decomposition in periodic components filled by closed orbits (two islands around centers and a cylinder in blue in the Figure) and minimal components two in the example)ads
%( Decomposition in periodic components filled by closed orbits (two islands around centers and a cylinder in blue in the Figure) and minimal components two in the example)
	
%	into  {periodic components} and {minimal components}. A  is a subsurface (possibly with boundary) on which all orbits are closed and periodic. These can be for example.  \emph{Minimal components} are subsurfaces (possibly with boundary) on which the flow is \emph{minimal} in the sense that all semi-infinite trajectories are dense 
% As proved independently by Maier \cite{Maier} and Zorich \cite{Zo:how}, $S$ can be decomposed The latter  are known as \emph{minimal components} of  $\varphi_\R$.   
%Also in the Nineties it was realized ( and independently Maier \cite{Maier}) that 

Minimal components of a locally Hamiltonian flow (and in particular minimal such flows, for which $S$ is in itself a minimal component), in suitably chosen coordinates, have the same orbits (up to \emph{time-reparametrization}, see \S\ref{sec:examples}) than  \emph{linear flows} discussed in  \S\ref{sec:linear} (see e.g.~\cite{Zo:how}). 

\subsection{Minimality and ergodicity}\label{sec:ergodicity}
To classify chaotic behaviour in locally Hamiltonian flows it is crucial to distinguish between two (complementary, up to measure zero) open sets (with respect to the topology described in \S~\ref{sec:generic}): in the first open set, which we will denote by $\mathscr{U}_{min}$, the typical flow is \emph{minimal} (in particular there are no centers and there is a unique minimal component). On the other open set that we will call $\mathscr{U}_{\neg min}$ there are periodic components (bounded by saddle loops homologous to zero), but the typical flow is still minimal when restricted to each complementary (minimal) component.  

 Let us remark that if the flow  $\varphi_\R$ given by a closed $1$-form $\eta$ has a \emph{saddle loop homologous to zero} (i.e.\ the saddle loop is a \emph{separating} curve on the surface), then the saddle loop is persistent under small pertubations (see Section~2.1 in \cite{Zo:how} or Lemma 2.4 in \cite{Rav:mix}). In particular, the set of locally Hamiltonian flows which have at least one saddle loop is open and gives the set denoted  $\mathscr{U}_{\neg min}$  above. The set $\mathscr{U}_{min}$ is given by the interior   (which one can show to be non-empty) of the complement of  $\mathscr{U}_{\neg min}$, i.e.\ the set of locally Hamiltonian flows without  saddle loops homologous to zero\footnote{Note that saddle loops non homologous to zero (and saddle connections) vanish after arbitrarily small perturbations and neither the set of 1-forms with saddle loops non homologous to zero (or saddle connections) nor its complement is open (see \cite{Rav:mix} for details).}.

The typical locally Hamiltonian flow (with respect to the measure defined in \S~\ref{sec:generic}) is $\mathscr{U}_{\neg min}$  is not only minimal, but also \emph{uniquely ergodic}. For a typical flow in $\mathscr{U}_{\neg min}$, the  restriction of the  flow on each minimal component is (uniquely) ergodic. Both results about minimality and ergodicity can be deduced from the classical results respectively by Keane and Masur and Veech (recalled in \S~\ref{sec:Teich}) respectively concerning of minimality and ergodicity of typical translation flows, by using that that the flow restricted to a minimal components is a time-change of a linear flow (see \ref{sec:linear})\footnote{One needs also to exploit that the two notions of typical, respectively for linear flows and locally Hamiltonian flows, interact well with each other: in particular, if a result holds for almost  linear flow in almost eveyr direction on almost every translation surface with respect to the Masur-Veech measure, one can show that it holds for a full measure set of locally Hamiltonian flows with respect to the Katok fundamental class.}.

\subsection{Mixing properties of locally Hamiltonian flows.}\label{sec:mixing} As mentioned earlier, finer chaotic properties such as (weak) mixing and spectral properties, crucially depend also on the \emph{speed} of motion along the orbits.  
The question of mixing in locally Hamiltonian flows was motivated by Arnold's conjecture in the $1990s$. In constrast with translation flows, which are never mixing (see \S~\ref{sec:Teich}), V. Arnold  in the $1990$s noticed a geometric phenomenon (explained in \S~\ref{sec:mixing}) which could produce mixing in locally Hamiltonian flows on the torus with one minimal component (those which we nowadays call~\emph{Arnold flows}, see \S~\ref{sec:minimal} and in particular footnote~\ref{Arnoldflowfootnote}). His intuition was proved to be correct shortly after by Sinai and Khanin \cite{SK:mix}, for a  full measure set (later improved by Kocergin \cite{Ko:mix})   of such flows on tori. The question of whether mixing is typical also for flows on higher genus %(when the singularities are non-degenerate)
 is much more delicate, and stayed open for two decades.

It turns out that mixing depends crucially on the type of singularities of the flow. When there are \emph{degenerate saddles} (i.e. multi saddles with $k \geq 6$ prongs, as in Fig.~\ref{multisaddle}), mixing had been proved already in the $1970s$ (by Kochergin in \cite{Ko:mix}) since in this case the saddles have a much stronger effect\footnote{As explained in \S~\ref{sec:shearing} mixing happens through shearing of transversal arcs and equidistribution of flow trajectories. In the case of degenerate saddles, the shearing effect is much faster and allows for simpler shearing estimates and faster mixing. Mixing in this case is indeed believed to have polynomial rates (see e.g.~a partial result in this direction in \cite{Fa1}), as opposed to the logarithmic (sub-polynomial) speed in the case of non-degenerate saddles (see \cite{Rav:mix} and the comments after Theorem~\ref{thm:mixing} in \S~\ref{sec:mixing}).}. In the  case of non-degenerate saddles (which, we recall, is \emph{generic} case, but much more delicate to treat), one has very different results in the open sets $\mathscr{U}_{min}$ and $\mathscr{U}_{\neg min}$ introduced in the previous \S~\ref{sec:ergodicity}. 
The full classification of mixing properties  has been a central part of my past research achievements \cite{Ul:mix, Ul:wea, Ul:abs}. The two following two results now give a complete picture:

\begin{theorem}[U' \cite{Ul:wea, Ul:abs}\label{thm:absence}\footnote{Absence of mixing for typical flows for any $g\geq 2$ was proved in \cite{Ul:abs}, while weak mixing is proved also for minimal components of locally Hamiltonian flow with simple saddles in \cite{Ul:wea}. Let us remark that a result in this direction for $g=1$ was already proved by Kocergin in \cite{Ko:abs} (see also \cite{Ko:abs2}) (in the language of special flows over rotations, which does not have a direct implication for locally Hamiltonian flows but suggested that the absence of mixing could hold also when rotations are replaced by IETs and hence in higher genus. Absence of mixing in the special case of $S$ with $g=2$ and locally Hamiltonian flows with two isomorphic simple saddles was shown by Scheglov \cite{Sch:abs}.}]
In $\mathscr{U}_{min}$, the typical locally Hamiltonian flow is \emph{weakly mixing, but it is \emph{not} mixing}. 
\end{theorem}

\begin{theorem}[U' \cite{Ul:mix}, Ravotti \cite{Rav:mix}\footnote{The result in this generality is proved in \cite{Rav:mix}. Mixing in a special (but crucial) case was proved in \cite{Ul:mix} in the language of special flows, more precisely for special flows over interval exchange maps under a roof with one asymmetric logarithmic singularity. The general case is the case of several logarithmic singularities with a global asymmetry condition (see e.g.~\cite{Rav:mix}).}]\label{thm:mixing}
In $\mathscr{U}_{\neg min}$, the restriction of the typical locally Hamiltonian flow $\varphi_\mathbb{R}$ on each of its minimal components is mixing. 
\end{theorem}

Let us remark even though the typical flow in  $\mathscr{U}_{min}$ is not mixing, there exists exceptional \emph{non-mixing} flow in this open set  as it was shown by J.~Chaika and A.~Wright \cite{CW}. Examples of \emph{mild mixing} (which is an intermediate property between weak mixing and mixing) were also built in \cite{KK} by A.~Kanigowski and  J.~Ku{\l a}ga-Przymus, but again are non typical (and one might conjecture that mild mixing is indeed non typical). 

Furthermore, there are also \emph{quantitative} results on the \emph{speed} of mixing (when there is mixing) which show that it happens (as expected in a parabolic flow) very slowly. More precisely, for a typical $\varphi_\R$  in $\mathscr{U}_{\neg min}$,  restricted to a a minimal component (which is mixing and hence display decay of correlations, refer to Section~\ref{sec:chaotic} for definitions), the \emph{speed} of decay of correlations (also sometimes called \emph{speed of mixing}) is \emph{sub-polynomial} (in accordance to what we expect for a slowly chaotic flow) and actually logarithmic,  namely for every pair $f,g $ of smooth observables there exists constants $c>0, \alpha>0 $ such that $|C_{f,g}(t)|\leq  c \log t ^\alpha$ (as shown by Ravotti in \cite{Rav:mix}).

\subsection{The role of shearing in slow mixing}\label{sec:shearing}
\begin{figure}[h!]
 \subfigure[shearing \label{shearing}]{
 \includegraphics[width=0.15\textwidth]{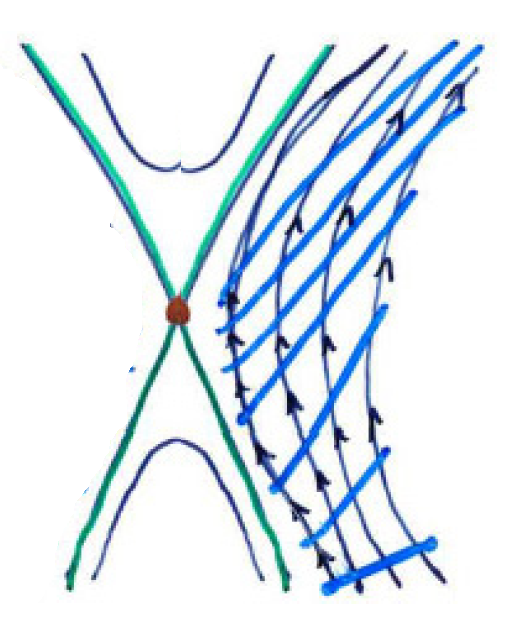}} \hspace{6mm}	
 \subfigure[compensation \label{compensation}]{
	\includegraphics[width=0.15\textwidth]{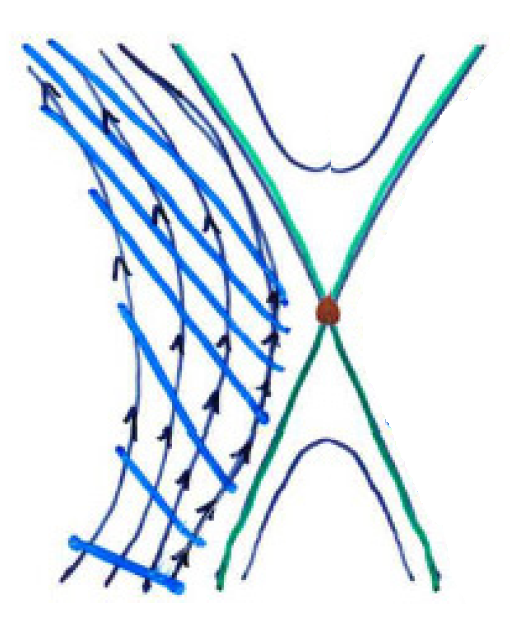}} \hspace{6.9mm}	
\subfigure[asymmetry \label{sheartrap}]{\includegraphics[width=0.2\textwidth]{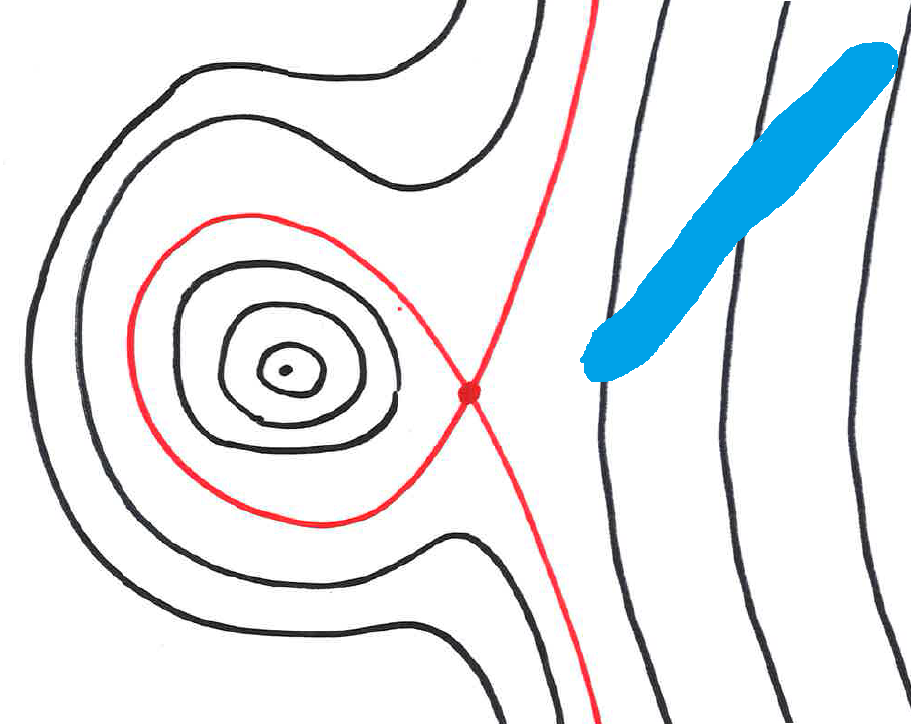}   }\hspace{6.9mm}	
\subfigure[wrapping \label{wrapped}]{
 \includegraphics[width=0.21\textwidth]{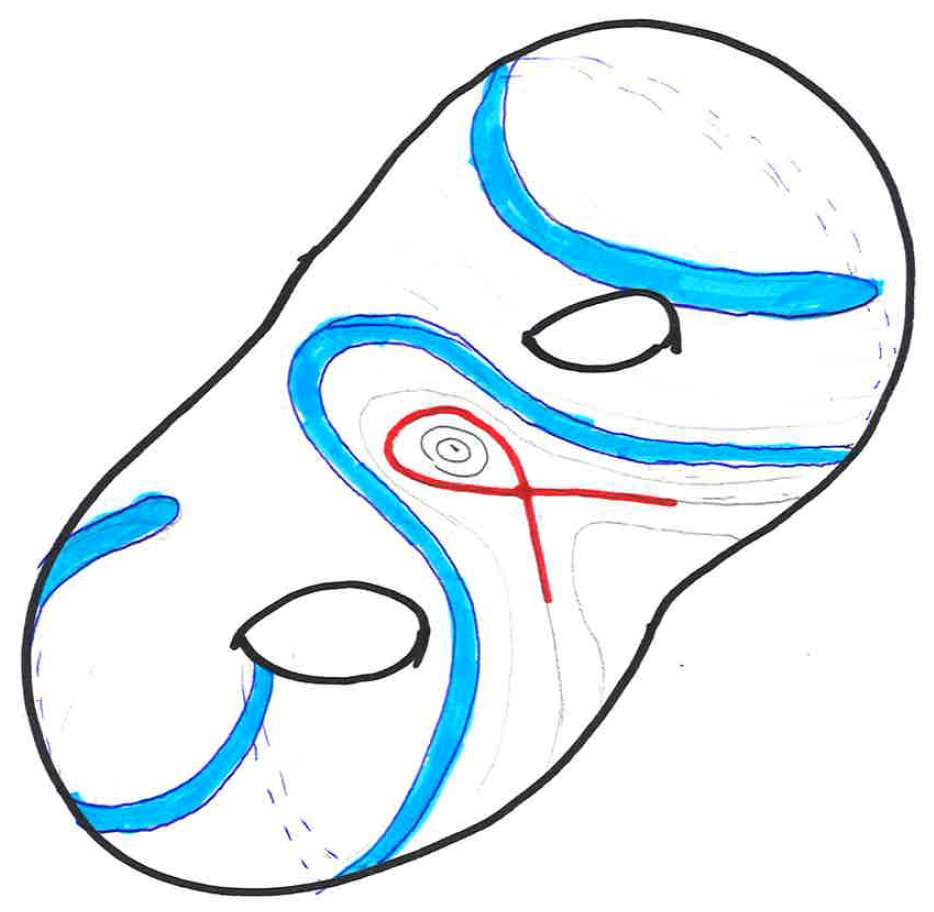}} 
 \caption{Shearing mechanism in locally Hamiltonian flows
.\label{shearing_mechanism}} 
\end{figure}

The parabolic nature of locally Hamiltonian flows is entirely due to the
presence of the saddles, which split nearby trajectories (as shown in Figure~\ref{fig:splitting}) and are responsible for the
slow divergence of nearby trajectories through a geometric phenomenon called \emph{shearing} (pictured in Figure~\ref{shearing_mechanism}). 
The butterfly effect in locally Hamiltonian flows happens indeed in a special way. 
 In presence of a Hamiltonian saddle, the closer a trajectory is to a saddle point, the more motion along the trajectory is slowed down. Thus, if we consider a small arc $\gamma$ transversal to the trajectories of the flow $\varphi_\R$, so that when flowing it,  $\varphi_t(\gamma)$ passes
nearby a saddle separatrix without hitting the saddle point (as shown in Figure~\ref{shearing_mechanism}), the different deceleration rates of points cause $\varphi_t(\gamma)$ to \emph{shear} in the direction of the 
flow (see Figure~\ref{shearing}). This description also evidences the slow nature of the butterfly effect in this case: points on nearby trajectories diverge from each other since, even though they travel on nearby trajectories, one travels faster than the other. Therefore the speed of divergence  is in this case only as fast as the speed of shearing. 

The shearing accumulated can  be later destroyed when the $\varphi_t(\gamma)$ passes near the other side of a saddle (see Figure~\ref{compensation}). The presence of a saddle loop, though, (as in Figure~\ref{sheartrap}) typically creates an asymmetry (this was the key intuition of Arnold that had motivated his conjecture on mixing) by producing stronger shearing on one side and hence,  in this case, the accumulation of shearing predominantly in one direction produces global shearing. 

%In the symmetric
%case, though, the shear which is gained when passing near one side of a singularity
%is lost when passing on a dierent side of other saddles.
% In the asymmetric case,the shearing happens prevalently in one direction, so it accumulates and after a
%long time, many small arcs transversal to the 
%ow become long sheared arcs almost
%in the 
%ow direction (see Figure 1) which wrap around the surface.

This geometric shearing phenomenon  is a crucial ingredient in the proofs of the mixing results in \S~\ref{sec:mixing} (in particular Theorem~\ref{thm:mixing}, but also to prove mixing in the exceptional examples in \cite{CW}) as it allows to deduce mixing from ergodicity (i.e.~equidistribution of flow trajectories, recall \S~\ref{sec:chaotic}). For large times $t\gg 0$, segments which do not hit the singularities will be so sheared along the flow, to be well approximated by long flow trajectories. Thus, for any given measurable set $A\subset X$, for every large $t$ one can cover an arbitrarily large proportion $A_t\subset A$ of  with a collection of short transversal segments $\{ \gamma_\alpha, \alpha \in \mathscr{A}_t \}$, \emph{each of which}, after time $t$, shadows a long trajectory of $\varphi_\R$, which is (close to)  equidistributed by ergodicity. One can hence show that each $\varphi_t(\gamma_\alpha)$ is also (close to) equidistributed.  Thus, since
$$
\varphi_t (A_t)\cap B = \cup_{\alpha\in\mathscr{A}_t} \left(\varphi_t(\gamma_\alpha)\cap B\right)
$$
by a Fubini argument, one can deduce equidistribution of $\varphi_t(A)$ (i.e.~mixing, see \S~\ref{sec:chaotic}) from equidistribution of each    $\varphi_t(\gamma_\alpha)$ (which follows as we said from shearing and unique ergodicity). Furthermore, the speed of mixing (or equivalently~the speed of decay of correleation) depends on the speed of shearing, which is \emph{slow} (namely subpolynomial in this case). 

\smallskip
This mechanims for \emph{mixing via shearing} seem to be a very common phenomen in slowly chaotic dynamics. 
The few of the early results on horocycle flows  (such as Marcus proof of mixing in \cite{Ma1} or Ratner's results, see \S~\ref{sec:Ratner}) exploit that small segments of geodesics curves, pushed by the horocycle flow, are \emph{sheared} in the horocycle direction\footnote{This shearing property follows from the commuting relations  between the horocycle flow  $h_\R$ and  the geodesic flow  $g_\R$, namely the relation $h_t g_s=g_sh_{e^{-2s}t}$ which holds for all $t,s \in \mathbb{R}$. The relation shows that if we push a small arc $\gamma = \{ g_s (x), s\in [0,\sigma]\}$ by $h_t$, the point $h_t (g_s(x))$ is \emph{aligned} in the geodesic direction not with $h_t(x)$, but with $h_{e^{-2s}t} (x)$ and hence that the pushed segments are \emph{sheared} in the flow direction.}. Furthermore,  since this is  essentially a \emph{geometric} mechanism for explaining mixing, this phenomenon \emph{persists} under perturbation and hence can be used also for time-changes (see \cite{Ma1} and \cite{FU}, where we prove quantitative mixing results and show polynomial estimates on the decay of correlations for smooth time-changes of the horocycle flow). 

A similar mechanism, namely shearing of segments of a suitable foliation (but with the difference that the direction of shearing is not global but depends on the segment considered) was also exploited in \cite{Fa:ana} to prove mixing in some (exceptional) elliptic flows\footnote{In \cite{Fa:ana}, Fayad shows that there exists mixing time-changes of linear flows on tori $\mathbb{T}^n=\mathbb{R}^n/\mathbb{Z}^n$ in dimension $n\geq 3$. The phenomenon is rare though, since it requires a highly Liouvillean rotation number.}  and in the context of nilflows: while nilflows are never mixing (see footnote \ref{nilflownevermixing}), in suitable classes of smooth time-changes one can implement this mechanism to prove mixing using shearing, see \cite{AFU, Rav:nil, AFRU}.
\smallskip

Finally, the complementary  results on \emph{absence}  of mixing (see Theorem~\ref{thm:absence}) involve showing \emph{absence of shearing}\footnote{To show absence of shearing, one needs to exploit that, when there are no saddle loops homologous to zero and hence no asymmetry which produces global shearing, the effect of shearing on two different sides of the same saddles compensate and cancels. This requires subtle estimates which hold for a full measure set of flows in $\mathscr{U}_{\neg min}$. In the exceptional mixing examples built in \cite{CW} sharing is still at the base of mixing, but is not produced by asymmetry of the singularities, but by an asymmetric equidistribution, so that trajectories, at different time scales, spend much more time on one side of a saddle than another)}. Indeed, a criterion for absence of mixing already formulated by Kocergin in \cite{Ko:abs} shows that (at least for typical) locally Hamiltonian flows mixing via shearing is essentially the \emph{only} possible way of achieving mixing.

\subsection{Beyond mixing, exploting shearing: Ratner's work}\label{sec:Ratner}
Whether one can deduce stronger and finer ergodic and spectral properties from shearing, in the context of flows with singularities, has been an open problem for decades, which has seen advances onlyl very recently (see \S\ref{sec:RatnerlocHam}). 
A great example of the fine and deep results on finer ergodic properties and rigidity phenomena  that one can obtain from shearing  is given by the celebrated works by Marina Ratner on the horocycle flow (and more generally unipotent flows in \emph{homogeneous dynamics}) \cite{Ra1, Ra2, Rat, Rat:coh}. Her work, and more in general 
the rigidity theory  for
%The ergodic and spectral theory of homogeneous parabolic flows (the 
 \emph{unipotent flows}, %(e.g. nilflows and unipotent flows)
% were laid in the 70s and 80s and
%  In particular, the theory of \emph{unipotent flows} 
 developed by Dani, Margulis and many others, has found breakthrough applications and has led to the solutions of important problems in number theory (such as the Oppenheim conjecture) and mathematical physics (such as the Bolztmann-grad limit for the Lorentz-gas). %; \emph{nilflows} have applications in number theory (such as distribution of fractional parts of polynomials and estimates of Weyl sums, see e.g.~\cite{Fla-Fo:nil}). 
%Unipotent flows are algebraic (homogeneous) flows  given by  \emph{uniponent} one-parameter subgroups  (in $ SL(n, %\mathbb{R})$ these are subgroups of matrices with all eigenvalues equal to $1$). Mixing for unipotent flows is well understood. Mixing unipotent flows  are  indeed mixing of all orders \cite{Ma}, and their spectrum is countable Lebesgue, in particular absolutely continuous (see \cite{Starkov} and the references therein). 
%One of the few classical results which applies to time-changes of unipotent flows is the proof by B.~Marcus  that all smooth time-changes of horocycle flows satisfying a mild differentiability conditions are mixing \cite{Ma}. 

\emph{Shearing} is at the at the  heart of Ratner's work and the above mentioned rigidity results.  A crucial ingredient in her work, indeed, is a technical  property introduced in \cite{Rat}  (that she calls property $H$), nowadays known  as the \emph{Ratner property} (see \cite{Thou}). It is this property, that Ratner  verified for horocycle flows, that is  used to deduce some of the main rigidity properties  of horocycle flows  (such as  joinings and measure rigidity).  % (originally Ratner called property $H$). 

 \begin{figure}[h!]
 \subfigure[Splitting (shift of $1-$time unit) and realignement.\label{RatnerShearing1}]{  \includegraphics[height=0.13\textheight]{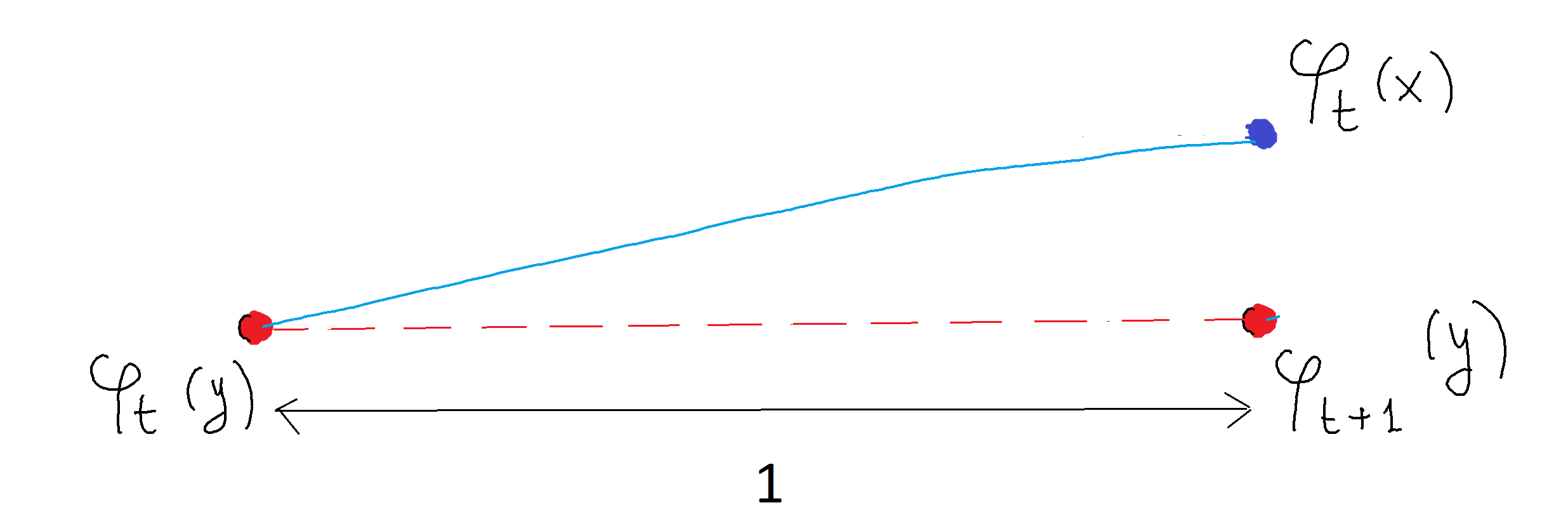}}\hspace{6mm}
 \subfigure[Realigned points stay close \label{RatnerShearing2}]{
  \includegraphics[height=0.13\textheight]{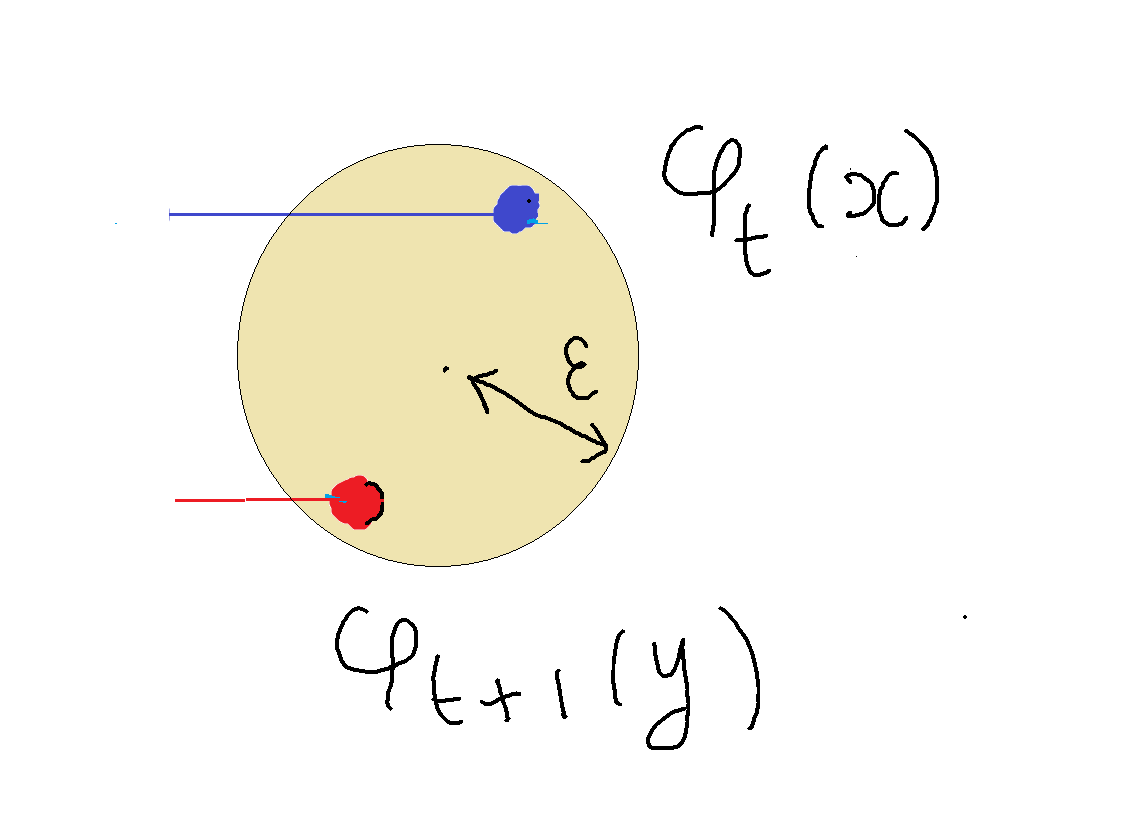} 	} 
%\subfigure[on $\mathbb{R}^2/\mathbb{Z}^2$ \label{Arnoldflowsquare}]{\includegraphics[width=0.33\textwidth]{LocHamSquare}   }	
 \caption{The Ratner property describing quantitative slow shearing. \label{RatnerShearing}}
\end{figure}

 %, whose formal is rather technical (see also ~\ref{sec:Ratner}) 
The \emph{Ratner property} encodes a quantitative property of \emph{controlled divergence} of nearby trajectories in the flow direction (illustrated in Figure~\ref{RatnerShearing}). 
%describes in a precise  quantitative way the shearing phenomenon, namely in which direction and how fast nearby trajectories diverge. 
% which is essentially a \emph{quantitative} description of t 
%Let $(X,d)$ be a $\sigma$-compact metric space, $\mathcal{B}$ the $\sigma$-algebra of Borel subsets of $X$ {and} $\mu$ a Borel probability measure on {$X$}.  Let $(T_t)_{t\in\mathbb{R}}$ be an ergodic flow acting on $(X,\mathcal{B},\mu)$.
Heuristically, it requires  that for \emph{most} pairs of nearby points $x,x'$, the orbits of $x,x'$ \emph{split} in the flow direction (say at time $t_1$) by a definite amount, called the \emph{shift} and then \emph{realign}, say by $\pm 1$ time-unit\footnote{One can more in general  consider a \emph{shift} $p$, which belongs to a fixed compact set $P$ so that now $\varphi_{t+p}(x)$ and the time-shifted orbit $\varphi_{t+p}(x')$ are close. The original Ratner property, where $P=\{ +1,-1\}$, is now sometimes called $2$-\emph{point} Ratner property, while the generalization to $P$ \emph{finite} first and to  $P$ \emph{compact}  
%finite first and compact later
later were defined by   by Fr\k{a}czek and Lema{\'{n}}czyk in \cite{FLRat, FLL} and called respectively \emph{finite} Ratner and \emph{weak} Ratner properties.}  so that now $\varphi_{t_1}(x)$ and the time-shifted orbit point $\varphi_{t_1\pm 1}(x')$ are close; then one requires the two orbits,  $\left(\varphi_{t_1+t}(x)\right)_{t\geq 0}$ and the time-shifted orbit $\left(\varphi_{t+1\pm 1+t}(x')\right)_{t\geq 0}$, to still stay close (see Figure~\ref{RatnerShearing}) for a fixed proportion $\kappa $ of the time $t_1$ it took to see the shift, namely for \emph{most} times $t \in [t_1 ,t_1+\kappa\, t_1]$.  
One can see that this type of phenomenon is possible only for parabolic systems, in which orbits of nearby points diverge with polynomial or subpolynomial speed.

\subsection{Searching for Ratner properties beyond unipotent flows}
%{\bf
%  %One might hope that the Ratner property, 
% We stress that the Ratner property holds by virtue of the polynomial divergence, or parabolic nature, of the horocycle flow.  
Since the Ratner property describes a form of divergence of nearby trajectories (or \emph{butterfly effect}) which is peculiar to parabolic flows, % (namely divergence happens via shearing, and with (sub)polynomial speed), 
it is reasonable to expect that some quantitative form of parabolic divergence similar to the Ratner property should hold  and be crucial in proving analogous rigidity properties for other classes of parabolic flows. Even more, since there is no formal definition for a system to be parabolic, one might even hope that the Ratner property could be taken as one of the characteristics making a system parabolic. 

The natural question hence arose  whether the Ratner property might hold   for \emph{smooth flows} on surfaces of higher genus. 
For a long time, though, there were \emph{no known examples} of systems with {the Ratner property} beyond horocycle flows and their (smooth) time changes. This changed drastically in the last decade. The first examples outside the homogeneous world were given by Fr\k{a}czek and Lema{\'{n}}czyk in \cite{FLRat, FLL, FL2d} (in the setting of special flows). The two authors could also show in  \cite{FLRat}  that a  \emph{variant}  of Ratner's property hold for some surface flows, more precisely in a class of flows on genus one tori known as \emph{von Neumann flows}\footnote{Von Neumann flows are called in this way since they first appeared in von Neumann's work \cite{vN}. In the special flows presentation they are flows over irrational rotations under a piecewise  linear roof with non zero sum of jumps.  For further recent results on von Neumann flows, see \cite{KS1, KS2, DK}.} (for non generic flows, corresponding to a  measure zero set of frequencies).  However, the flows in \cite{FLRat} are not (globally) smooth. % These are built from  genus one surfaces with boundary
 %For other examples with this type of Ratner property, see also  \cite{FLL, FL2d}. 

%another6 variation of the Ratner property was also showned by the same authors to hold for flows on $3$-dimensional tori, see \cite{FL2}). 
%K.\ Fr\c{a}czek and the second author showed in \cite{Fr-Le} that there exists a class of flows on $\T^2$ (smooth flows with one singular point) for which a variant of {the R-property} holds. 

 \begin{figure}[h!]
 \subfigure{  \includegraphics[height=0.11\textheight]{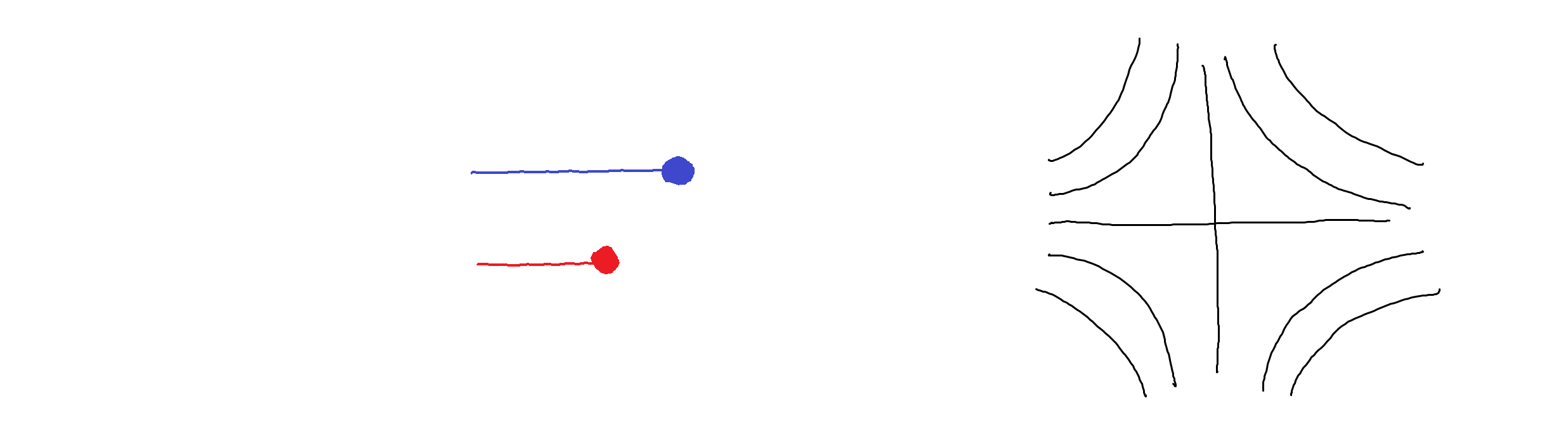}}\hspace{6mm}
 \subfigure{
  \includegraphics[height=0.11\textheight]{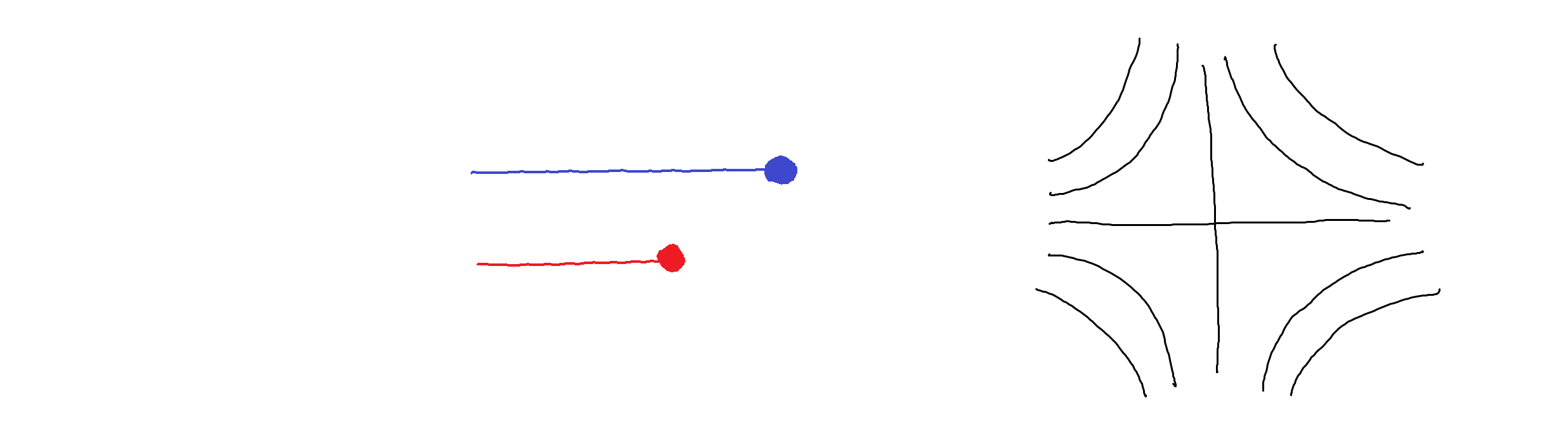} 	} 
%\subfigure[on $\mathbb{R}^2/\mathbb{Z}^2$ \label{Arnoldflowsquare}]{\includegraphics[width=0.33\textwidth]{LocHamSquare}   }	
 \caption{Splitting of trajectories of a  locally Hamiltonian flow near a saddle. \label{fig:splitting}}
\end{figure}

The difficulty in treating smooth flows on higher genus surfaces is given by the presence of singularities (which are unavoidable when $g\geq 2$, see \S~\ref{sec:locHam}), which introduce discontinuities and destroy the slow form of divergence {a la Ratner}: essentially, as soon as two nearby trajectories are separated by hitting a saddle (see Fig.~\ref{fig:splitting}), one drastically looses control of the divergence.  The Ratner property in its classical form (as well as the weaker versions defined in \cite{FLRat, FLL}) is  expected to fail for of smooth area-preserving flows with  non-degenerate  fixed points\footnote{The failure of the classical Ratner property was formally proved in a special case in ~\cite{FK} (for a  class of Kochergin flows, i.e.~special flows with power singularities over rotations, see Theorem 1 and the Appendix B in ~\cite{FK}) and this result gives reasons to believe that,  for similar reasons, the classical Ratner property should indeed always  fail in presence of singularities.}.

% Unfortunately, when parabolic flows have singularities, as it happens for smooth area-preserving flows with fixed points, the
% The heuristic problem for Arnold flows and more generally smooth area-preserving flows to enjoy the Ratner property (or its weaker versions) 
% is that Ratner-like properties require a (polynomially) controlled way of divergence of orbits  of nearby points, but if the orbits of two nearby points get too close to a singularity, their distance explodes in an uncontrolled manner (see e.g. the Appendix of \cite{FK}).   

\subsection{The Switchable Ratner property in locally Hamiltonian flows.}\label{sec:RatnerlocHam}
The possibility of pushing our understanding of smooth flows on surfaces, using techniques loosely inspired by Ratner's work, emerged only recently, in virtue of the recent developments in the field. 
 %The results in  \cite{KKU} provide the first example in which (a variant of)  the Ratner property was proved outside of the homogeneous world, for a class of \emph{typical} smooth surface flows.   
 A key breakthrough was achieved recently by B.~Fayad and A.~Kanigowski, who, in \cite{FK}, introduced a new modification of the Ratner property, the so called \emph{Switchable Ratner} property (or \emph{SR-property}). According to this variation, it is sufficient to see the Ratner divergence of orbits for most pairs of initial conditions $(x,y)$  \emph{either}  in the future (for $t>0$) \emph{or} in the past (for $t<0$), depending on the pair of initial points. Thus, if one pair of nearby trajectories is separated by hitting a singularity (as shown in Fig.~\ref{fig:splitting}), and hence their distance explodes in an uncontrolled manner, one can still hope to be able to prove the Ratner slow form of divergence when flowing backward in time.% (see Figure~\ref{switchableRatner}). 

Let us remark that  above mentioned variations  Ratner property (thus in particular also the switchable Ratner property) were defined in order to have the same strong dynamical consequences of the original Ratner property. In particular, all variants of the Ratner property, as the original Ratner property does, imply a \emph{rigidity}-type result on  \emph{joinings}\footnote{The notion of \emph{joining} plays a key role in ergodic theory. Assume that $\varphi_\R,\phi_\R$ are flows on probability standard Borel spaces $(X, \mu)$ and $(Y,\nu)$, respectively. A \emph{joining} between  $\varphi_\R$ and $\phi_\R$  is a  $\varphi_\R\times \phi_\R$ -invariant (for each $t\in\R$) probability measure on $X\times Y$ with the projections $\mu$ and $\nu$, respectively.\label{footnote:joiningdef} 
By $J(\varphi_\R, \phi_\R)$ we denote the set of {\rm joinings} between the flows $\varphi_\R$ and $\phi_\R$.  A trivial joining always exists and is given by the product measure $\mu\times \nu$.}
%Following Furstenberg \cite{Fu}, we say that $T_t$ and $S_t$ are {\em disjoint}, and we write $T_t\perp S_t$, if  $\mu\otimes\nu$ is the only member of $J(T_t,S_t)$.}
% Note that if $T_t\perp S_t$ then at least one of the two flows must be ergodic. If both flows are ergodic, then by $J^e(T_t,S_t)$ we denote the subset of {\em ergodic joinings},  that is, of those $\rho\in J(T_t,S_t)$ which make the flow $(T_t\times S_t)_{t\in\R}\subset{\rm Aut}(X\times Y,\mathcal{B}\otimes\mathcal{C},\rho)$ ergodic (this set is always nonempty, as the ergodic decomposition of any joining yields a.e.\ ergodic component also a joining). Note also that if $W$ yields an isomorphism of $T_t$ and $S_t$ then the formula
%$$
%\rho(B\times C):=\mu(B\cap W^{-1}C),\;\;B\in\mathcal{B},C\in\mathcal{C},$$
%determines a member of $J(T_t,S_t)$ which is ergodic if $T_t$ (hence $S_t$) is ergodic.  Note also that, for each $r\neq0$, we have $J(T_t,S_t)=J(T_{rt},S_{rt})$ with the equality $J^e(T_t,S_t)=J^e(T_{rt},S_{rt})$ whenever the flows $T_t$ and $S_t$ are ergodic.}) 
results (by restricting the type of
\emph{self-joinings} that the flow can have\footnote{The Ratner property and its variations in particular imply a property known as \emph{finite extension of joinings}:  any  non-trivial ergodic joining $\rho $ is a finite extension  (finite union of graphs).\label{footnote:joiningrigidity}  Furthermore, Ratner properties impose some restrictions not only on the set of self-joinings but also on the set of its joinings with another (ergodic) flow. In particular, we can ask whether for two systems {sharing the same  Ratner property } it is possible to \emph{classify} joinings between them. The first result in this direction can be found { in Ratner's work~\cite{Rat:coh}, where she shows} that two flows $\bar{h}_\R$ and $\tilde{h}_\R$ given by two smooth time changes of a horocycle flow $h_\R$  are disjoint, i.e.\ the only joining between them is the  product measure, whenever the cocycles corresponding to the time changes are not {\it jointly cohomologous} (see~\cite{Rat:coh} for definitions).}) and 
allow to enhance mixing properties (see for example Corollary~\ref{cor:mixingallorders}).

B.~Fayad and A.~Kanigowski could prove in \cite{FK} that this variation of the Ratner property holds for some smooth surface flows in genus one, more precisely for typical Arnold flows (see~\ref{sec:minimal}, Figure~\ref{Arnoldtorus})  as well as for (a measure zero class of) torus flows with one degenerate (or fake) singularity  (sometimes known as \emph{Kocergin flows}). Let us recall that in higher genus ($g\geq 2$) it is important to distinguish between the two open sets $\mathscr{U}_{min}$ and $\mathscr{U}_{\neg min}$ (see \S~\ref{sec:ergodicity}) of locally Hamiltonian flows with non-degenerate singularities. 
 In  \cite{KK}, the SR-property was proved for some (measure zero set of)  minimal smooth flows \footnote{The flows considered in  \cite{KK} are special flows over IETs of \emph{bounded type}, under a roof with symmetric logarithmic singularities. Bounded type Diophantine conditions on IETs extend the notion of bounded type (also called \emph{constant type}) rotation numbers and, as among rotations, are measure zero (but full Hausdorff dimension) conditions.} in $\mathscr{U}_{min}$. It is likely that to prove a form of Ratner properties for other  flows  (hopefully a full measure set) in $\mathscr{U}_{min}$ will require introducing yet another variant of the Ratner property, one which could take into consideration \emph{average} shearing and thus will require new ideas.

The result in \cite{KK}, on the other hand, shows that the switchable Ratner property holds for (the minimal component of) the typical (Arnold) flow in $\mathscr{U}_{\neg min}$ when $g=1$ and the flow has only one simple saddl (and center).  In joint work with A. Kanigowski and J.~Ku{\l a}ga-Przymus \cite{KKU}, we could prove that the \emph{switchable} version of the Ratner property  is \emph{typical} among mixing (components of) locally Hamiltonian flows in    $\mathscr{U}_{\neg min}$  for  \emph{any genus} $g\geq 1$ (thus extending to more singularities\footnote{For $g=1$, B.~Fayad and A.~Kanigowski show in \cite{FK} that the switchable Ratner property holds for (the minimal component of) a full measure set of Arnold flows only when there is a unique saddle; they also consider the case of more saddles, but then require a condition which has measure zero. The main result in  \cite{KKU} (which can be expressed in the language of special flows over IETs under roofs with asymmetric logarithmic singularities) on the other hand  gives a full measure condition not only the higher genus case, but also the case of genus one and more saddles.} and generalizing to higher genus $g\geq 2$ the result by \cite{FK}):

\begin{theorem}[Kanigowski, Ku\l aga-Przymus and U', \cite{KKU}]
For any $g\geq 1$, a typical locally Hamiltonian flow $\varphi_\R$ in $\mathscr{U}_{\neg min}$, restricted on any of its mininal component, has the switchable Ratner form of shearing. 
\end{theorem}
This result hence imply a rigidity type result for the classification of joings (see footnote \ref{footnote:joiningrigidity}) and in particular allowed us to upgrade mixing to  a stronger property, namely mixing of all orders (see \S~\ref{sec:chaotic} and \eqref{def:multiplemixing} for the definition).  

\begin{corollary}[KKU]\label{cor:mixingallorders}
For any $g\geq 1$, the restriction of a typical locally Hamiltonian flow $\varphi_\R$ in $\mathscr{U}_{\neg min}$ on any of its miminal components is mixing of all orders.
\end{corollary}

Thus, the Corollary show that Rohlin's conjecture (see Section \ref{sec:chaotic}) holds for these class of smooth flows.

Further recent works also show that Ratner properties also hold for other classes of slowly chaotic flows. For example the Switchable Ratner property holds for a class of time changes of constant type {\it Heisenberg nilflows}, see the recent work \cite{FK} by Forni and Kanigowski.

\subsection{{\bf Disjointness of rescalings.}} \label{sec:disjoint}
Advances in our understanding of disjointness properties became possible building on the switchable Ratner property \cite{FU, FFK, BK}. The notion of \emph{disjointness}\footnote{Two measure preserving flows $\varphi_\R$ and $\phi_\R$ are called \emph{disjoint} (in the sense of Fursterberg) if their only common ergodic joining is the trivial joining (i.e.~the product joining).} was introduced in the $1970s$ by H.~Furstenberg (see in particular \cite{FuR}); two disjoint flows are in particular \emph{not isomorphic}\footnote{We say that two flows $\varphi_\R$ on $(X,\mu)$ and $\phi_\R$ on $(Y,\nu)$ and are \emph{isomorphic} as measure preserving flows if there exists a isomorphism $\Phi:X\to Y$, i.e.~a one-to-one measureable map which respects measure zero sets and commutes with the dynamics, i.e.~$\Phi\circ \phi_t(x)=  \psi_t(\Phi(x))$ for $\mu$-almost every $x\in X$. If two flows are isomorphic,  the isomorphism $\Psi$ yields a non-trivial joining, so that the two flows cannot be disjoint.}.

A disjointness property which has received a lot of attention recently (in particular as a tool in connection with Sarnak's conjecture on Moebius orthogonality, see below) is \emph{disjointness of rescalings}.  
Given a real number $\kappa>0$, by the $\kappa$-\emph{rescaling} of  $\varphi_\R$ we simply mean the flow $\varphi^\kappa_\R:=(\varphi_{\kappa t})_{t\in \R}$ (in which the time is rescaled by the factor $\kappa$)\footnote{Notice that when $p$ is an integer, the time-one map of $p$-rescaling coincides with the $p$-\emph{power} of the time-one map $R_1$ of the flow, so considering rescalings is an analogous operation to considering the powers of a given transformation.}. Thus, a rescaling is a special type of time-reparametrization of a flow, given by a \emph{linear} time-change. We say that  $\varphi_\R$  has \emph{disjoint rescalings} if for all (or all but finitely many) $p,q>0$, the rescalings $\varphi_\R^p$ and $\varphi_\R^q$, where $p,q>0$ and $p\neq q$, are disjoint (in the sense of Furstenberg). Disjointness of rescalings has played a key role in proving some of the first  instances of Sarnak's conjecture \cite{Sa} of orthogonality of the Moebius function\footnote{Let us recall that the \emph{Moebius function} $\mu:\mathbb{N}\to \mathbb{Z}$ is a multiplicative function defined on a prime $p$ by $\mu(p):=({-1})^k$ if $p=p_1\cdots p_k$ with $p_i$ \emph{istinct primes}. Sarnak conjectured (see \cite{Sa}) that $\mu $ is \emph{orthogonal} to entropy zero deterministic sequences, i.e.~if $f:X\to X$ is a topological dynamical systems (i.e.~$X$ is a topological space and $f$ is continous) then, for every $x\in X$, $\sum_{n=1}^N f(T^n x ) \mu(n) = o(N)$ as $N$ tends to infinity. The survey \cite{surveyMoebius} gives an introduction to the conjectures and the progress made so far.} in number theory with entropy zero dynamical systems (as a tool to prove the conjecture via the so called Katai orthogonality criterion, see for example \cite{Bo-Sa-Zi, FFnote} and more in general the survey \cite{surveyMoebius}). 
% \cite{Katai} and \cite{}
%One more motivation for studying joinings rather then just self-similarity of different rescalings of
%a given flow {Tt}t∈R has its justification in , which has
%recently been proven to be an important tool for studying problems around Sarnak’s conjecture on
%Möbius disjointness (see e.g. [9]).  

\smallskip
In recent joint work with A.~Kanigowski and M.~Lema{\'{n}}czyk \cite{KLU}, we introduced a new tool to study disjointness phenomena for smooth surface flows, namely a  disjointness criterion based on  the \emph{switchable} Ratner property.  The criterion was devised and formulated so that it can be applied  to prove disjointness of two flows which both have the switchable Ratner property, so that in both flows one can observe a controlled form of divergence of nearby trajectories (for example polynomial divergence), but  the  speed of divergence for the two flows is different
%. More precisely, we have two flows for which we can observe divergence of nearby trajectories in the flow direction in a controlled way (say with polynomial speed) but the speeds are different
(for example for one flow is it  linear, in the other quadratic).\footnote{ 
The key \emph{shearing phenomenon} exploited in the criterion is that, for pairs of two nearby points in the first system  and two nearby points in  the second, after some time (depending on both pairs of points) we will see a \emph{relative divergence}, i.~e.~in one pair we will see a realignement with some shift $C_1>0$ in the flow direction, while in the other we will see a realignement with a shift $C_2>0$ with $C_1\neq C_2$. This explains how the relative shearing appears in this situation.}
%Let us remark that if both flows have only the \emph{switchable} Ratner property, one might have pairs of points for which one see the Ratner form of shearing only in the past for one flow, and only in the future for the other. Therefore, there is no hope to implement the heuristic {described} above on the \emph{whole} (or large parts) of the space.
%An essential feature of the criterion is that for one of the two flows one can consider only pairs of points  in \emph{small} parts of space (whose  measure is  bounded {from} below, but not necessarily close to one), on which one sees controlled shearing \emph{both} forward \emph{and} backward. Thus, one can couple these pairs with pairs in the other flow which have the Ratner property in the future \emph{or} the past to get the \emph{relative divergence} explained above. We will come back to this point after the statement of the criterion in  {Section}~\ref{s:Sec3}.
%In order to give the precise formulation of the criterion, we now first need to give the definition of \emph{almost linear} reparametrization. % (which will be used in the criterion to describe the realignement of the trajectories with the smaller.
 Exploiting this criterium, we were able to show that disjoitness of rescalings is typycal among  Arnold flows (see \ref{sec:minimal}, Figure~\ref{Arnoldtorus}):

\begin{theorem}[Kanigowski-Lemanczyk-U', \cite{KLU}] \label{thm:KLU}
A {\emph{typical}} {Arnold flow}  has {disjoint rescalings}. In particular, if $\varphi_\R$ is the restriction of the unique minimal component of an Arnold flow, there exists only two values\footnote{The two values are related to the \emph{asymmetry} of the saddles: in the special flow representation, the Arnold flow can be written as a flow over a rotation $R_\alpha :[0,1)\to [0,1)$ (given by $R_{\alpha}(x)=x+\alpha\mod 1$) under a roof with logarithmic singularities given by the roof function $r(x)=C_0 \log |x| + C_1 \log (1-x)+h(x)$ where $h$ is smooth on $[0,1]$; then $q = C_0/C_1$.} of the form  $q,1/q$ such that $\varphi_\R$ and $\varphi_\R^p$ are disjoint for any positive $p\notin \{1,q,1/q\}$.  
\end{theorem}
As a Corollary, Sarnak's Moebius disjointness conjecture holds for these flows (see \cite{KLU} for details).

\smallskip
We believe disjoitness of rescalings should also hold for typical locally Hamiltonian flows  in higher genus, but this is currently an open problem. Preliminary work seems to indicate that, despite some technical additional difficulties, the techniques used to prove Theorem~\ref{thm:KLU} should allow to prove disjointness of rescalings for all minimal components of typical flows in the open set $\mathscr{U}_{\neg min}$ (refer to \S~\ref{sec:ergodic} for the definition of the open set $\mathscr{U}_{\neg min}$ ). The recent work  \cite{BK} by Berk and Kanigowski, even though it does not apply to surface flows directly, gives a good indication that disjointness of rescalings could also hold for typical flows (under a suitable full measure Diophantine-type condition) in the complementary set $\mathscr{U}_{\neg min}$. 

It is natural to ask whether disjoitness of rescalings could actually be a widespread feature of slowly chaotic systems. 
%We believe disjoitness of rescalings should also hold in higher genus, i.e.~for typical locally Hamiltonina flows in the open set $\mathscr{U}_{\neg min}$ (defined in \S~\ref{sec:minimal}), but this is currently still open. There are also reasons
The new disjointness criterion is also used in \cite{KLU} to prove disjointness of rescalings for (a class of) smooth (non-trivial) time-changes of the horocycle flow (see also \cite{FFnote}, where the result is proved for a more general class of time changes with different methods), answering  in particular a question of  Marina Ratner. Notice  that {two different} rescalings of the (classical, non time-changed) horocycle flow $h_\R$ are \emph{never} disjoint\footnote{Indeed,
if $(g_s)$ denotes the \emph{geodesic flow}, 
the renormalization equation $h_t g_s=g_sh_{e^{-2s}t}$ for all  $t,s \in \mathbb{R}$, 
yields that, for any positive $p\neq q$, the flows $h_\R^p$ and $h^q_\R$ are conjugated by $g_s$ with $s=-\frac{\log(q/p)}{2}$ (and hence are not disjoint).}. Thus, it seems that (non trivial) time-changes of the horocycle flow are in some sense better behaved and display chaotic features that the horocycle flow itself lacks,  due to its homogeneous and self-similar nature. Perhaps unfortunately, the most studied and best understood  model of a parabolic flow, the horocycle flow, may have not provided the most significant example in terms of generic chaotic properties. Hence, the importance of better understanding new and larger classes of slowly chaotic systems and their typical chaotic features.

% we explain the criterion and some of the applications drawn from it in \cite{KLU}. 
The new criterion for disjointness introduced in \cite{KLU} has already proved usesul in different contexts, see for example the recent works \cite{FoK, DK} where it is applied  to study disjointess phenomena respectively  for Heisenberg nilflows in \cite{FoK}, for von Neumann flows in genus one (as well as other special flows over IETs) in \cite{DK}. %and in \cite{BK} to prove disjointness of rescalings in some special flows over rotations.  
Finally,  the disjointness criterion is used in \cite{KLU} also to show that a typical Arnold flow is disjoint from  any smooth time change of the horocycle flow (and in particular from the classical horocycle flow itself), thus showing that  
these  two classes of parabolic flows are truly distinct.

%We conclude by saying that the (parabolic) disjointness criterion (Theorem 3.1) seems
%to be a general tool in studying disjointness of systems with a Ratner property. A variant
%of this criterion is used in a recent work of G. Forni and the first author [14], where
%disjointness properties of time changes of Heisenberg nilflows are studied and also in recent
%works of P. Berk and the first author [4] and C. Dong and the first author [8] for getting
%disjointness of some classes of interval exchange transformations (IET’s) and translation
%flows

\subsection{Spectral theory of locally Hamiltonian flows.}\label{sec:spectral} 
The last aspect we want to discuss is the \emph{spectral theory} of smooth area-preserving flows (introduced in \S~\ref{sec:chaotic}).  
The spectral properties (and in particular what is the spectral type, see \S\ref{sec:spectral} for definitions) of locally Hamiltonian flows is a natural question, which has been lingering for decades (see e.g.~\cite[Section 6]{KT} and \cite{L}).  
While the classification of mixing properties of locally Hamiltoninan flows is essentially complete (as summarized in \S~\ref{sec:mixing}), very little is known on the spectral properties of these flows beyond  the case of genus one. % Results on the spectrum of the operator, though, are very rare.

 One of the few results   in the literature concerning  on the nature of the spectrum of some area-preserving surface flows was proved by K. Fr\k{a}czek and M. Lema{\'{n}}czyk  in \cite{FL1}, who showed singularity of the spectrum  for \emph{Blokhin examples}\footnote{ In the setting of special flows,  Fr\k{a}czek and Lema\'nczyk in \cite{FL03}, study  special flows over rotations with single symmetric logarithmic singularity.  In \cite[Theorem 12]{FL03} it is shown that (for a full measure set of rotation numbers) such special flows are \emph{spectrally disjoint} from all mixing flows, from which it follows in particular that the spectrum is purely singular.} (which were the first examples of ergodic flows on surfaces, see \cite{Bl}). This gives examples of locally Hamiltonian flows on surfaces of any genus $\geq 1$ with  singular continuous spectrum  (see \cite[Theorem 1]{FL03}), but these are essentially built glueing genus one flows and thus, they are highly non typical.  
%(as special flows over rotations under a roof with symmetric logarithmic singularities, see \S\ref{sec:ownM}  for definitions). % which is closely related to non-degenerate Hamiltonian flows in genus one . fexa
 
\smallskip
It turns out that the geometric approach to (quantitative) mixing through shearing can sometimes be pushed to provide also spectral information on parabolic flows. For example,% quantitative  
% Pushing this geometric approach, % and combining it with quantitative results, 
in joint work with Forni \cite{FU} we were able to show that from \emph{shearing estimates} on can also access information about the \emph{spectrum} of time-changes of the horocycle flow, in particular settling in particular a conjecture of Katok-Thouvenot  (see also   \cite{Tie} where Tiedra de Aldecoa gave simultaneously a different proof by operator methods. {Further developments based on our approach were achieved also for unipotent flows by Lucia Simonelli  \cite{Si}. 

% was achieved by Fayad, Forni and Kanigowski in \cite{FFK}, who showed that a class of smooth flows on surfaces of genus one (which can also be represented as special flows over rotations, see \S\ref{sec:history}) has countable Lebesgue spectrum. These flows display a strong form of \emph{shearing} of nearby trajectories and
% were proved to be mixing by Kochergin  in the 70's, \cite{Ko:mix}.

A  breakthrough on the spectral side  was achieved recently in \cite{FFK} for a special class of smooth area preserving flows with \emph{degenerate} singularities \emph{when the genus of the underlying surface is one}, sometimes known as \emph{Kochergin flows} (since Kochergin \cite{Ko:mix} proved  their mixing, for any $g\geq 1$). These are minimal flows on the torus with one stopping point (see Fig.~\ref{Kocergintorus}), also called  \emph{fake} singularity (this point can be seen as a degenerate fixed point with only $k=2$ prongs).    Taking as a starting point the idea used by Forni an myself in \cite{FU}  to prove absolute continuty of the spectrum for time changes of horocycle flows, Forni, Fayad and Kanigowski, proved that, if the  \emph{degenerate} singularity is sufficiently strong\footnote{Kocergin flows admit a special flow representation where the roof has power-type singularities. In genus one, for a flow with one degenerate singularity, one has a special flow over a rotation, under the roof $r(x)=c_0/x^{\alpha}+c_1/(1-x)^\alpha$ for some $0<\alpha<1$. The assumption in \cite{FFK} is that $\alpha$ is  sufficiently close to $1$.}, the spectrum is \emph{absolutely continuous} (and actually \emph{Lebesgue}).
  
  \begin{theorem}[Fayad, Forni, Kanigowski, \cite{FFK}] \label{ac}
A locally Hamiltonian flow in \emph{genus one} with only one  sufficiently strong degenerate singularity as fixed point has {\emph{countable Lebesgue}}  {spectrum}.
\end{theorem}
Countable Lebesgue spectrum is a strong spectral result, which implies in particular that the spectrum is absolutely continuous (see Section~\ref{sec:chaotic} for the definition).   The result provides the first example of such a strong chaotic property in an entropy zero and low dimensional smooth system.  We remark that stopping points (and more in general degenerate fixed points)  are known to produce shearing and hence mixing \cite{Ko:mix} (at rates which are expected to be polynomial, see e.g.~\cite{Fa1}). The absolute continuity of the spectrum is essentially\footnote{There are actually technical issues in controlling the part of space where there is not enough shearing.} based on a decay of correlations which is square-summable. 

 \begin{figure}[h!]
  \subfigure[Kochergin flow (a.c.~spectrum)\label{Kocergintorus}]{  \includegraphics[width=0.3\textwidth]{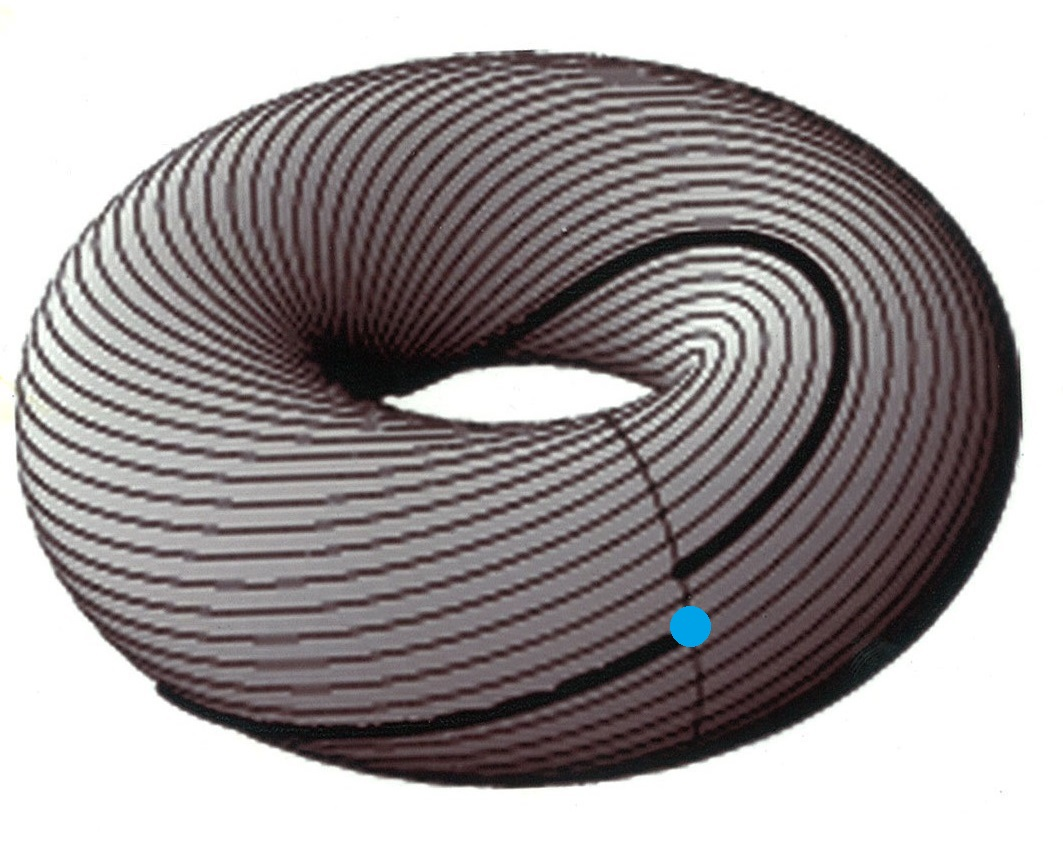}} \hspace{6mm} \subfigure[$g=2$, two isomorphic simple saddles (singular)\label{onlysaddles}]{ \includegraphics[width=0.5\textwidth]{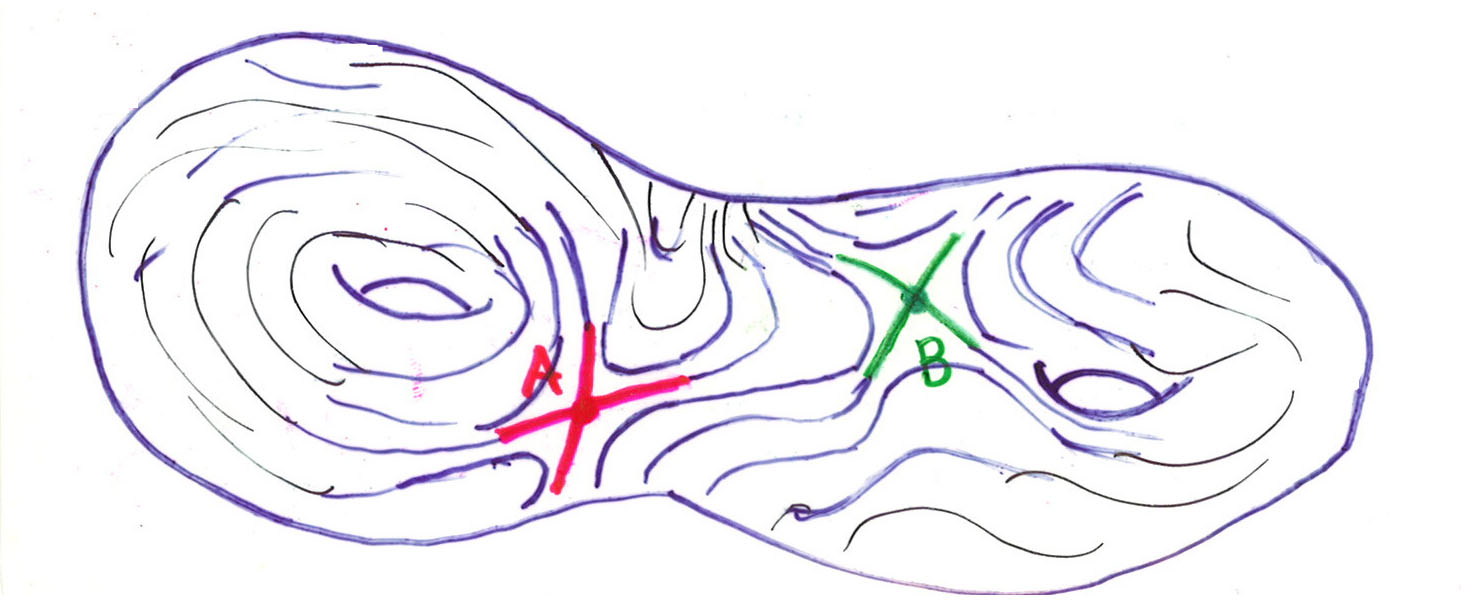}}\hspace{9mm}

{\caption{The locally Hamiltonian flows with absolutely continuous and singular spectra respectively in Theorems\ref{ac} and \ref{singular}.\label{saddles}}}
\end{figure}

A recent spectral breakthrough, which goes in the opposite direction, concerns the nature of  the spectrum of locally Hamiltonian flows on genus two surfaces, and, to the best of our knowledge, is the first general spectral result  for surfaces of higher genus, namely $g \geq 2$.   
%The first result, to the best of our knowledge, for higher genus surfaces, i.e. for $g\geq 2$, is the following result for genus two (see Figure~\ref{onlysaddles}), which goes in a completely opposite direction.
\begin{theorem}[Chaika-Fraczek-Kanigowski-U', \cite{CFKU}] \label{singular}
A {\it typical} locally Hamiltonian flow on a \emph{genus two} surface with two isomorphic {\it simple saddles} has {\emph{purely singular}}  {spectrum}.
\end{theorem}
%references\cite{FL1}
This result in genus two was inspired by the singularity result proved by  Fr{\c a}czek and M. Lema{\'{n}}czyk (for special flows over rotations) in \cite{FL03}. Their result indeed shows that, when one can prove absence of mixing and  some form of (partial) rigidity, it might be possible to deduce singularity of the spectrum.  
%A previous result which inspired some of the techniques (in particular the singularity criterion) used to prove Theorem~\ref{CFKU} was proved by  Fr{\c a}czek and M. Lema{\'{n}}czyk (for special flows over rotations), but, (to the best of our knowledge) singularity not known before for any \emph{smooth} surface flow in higher genus. 
%Even though singularity was   Fr{\c a}czek and M. Lema{\'{n}}czyk,
%\cite{FL05}, \cite{FL03} %%FL03 FL05
Theorem~\ref{singular} %can be also seen as a generalization o \cite{FL03,Fr-Le2}) 
 strenghthens one of the early results on absence of mixing, i.e.~the absence of mixing for typical flows in  the same class ($g=2$, two isomorphic saddles) proved by Scheglov \cite{Sch} (which is a special case of Theorem~\ref{thm:absence}). As in  \cite{Sch}, the assumptions are crucial since the underlying surface has an inner symmetry\footnote{More precisely, the linear flow of which the locally Hamiltonian flow is a time-change is a flow on a translation surface $S$ which admits an \emph{hyperelliptic involution}, i.e.~an affine automorphism $\Phi: X\to X$ which is an involution, i.e.~$\Phi^2=Id$.} which plays a crucial role in the proof. Nevertheless, we believe it should be possible to extend the result on higher genus exploiting the same singularity criterion used in \cite{CFKU} (which is an extension of the  criterion used in \cite{FL03} as well as \cite{FL05}), coupled with the delicate estimates for absence of shearing proved in \cite{Ul:abs} to prove Theorem~\ref{thm:absence}. 

\smallskip
The nature of the spectrum for other classes of locally Hamiltonian flows is unknown. It might be conjectured, in view of Theorem~\ref{ac} for Kocergin flows, that also in higher genus, in presence of sufficiently strong \emph{degenerate} singular points, the spectrum is also absolutely continuous (and even countable Lebesgue), essentially thanks to a strong quantitative control of decay of correlations.  It is not clear what to expect when the degenerate singularity is not sufficiently strong\footnote{As already remarked in a footnote, sufficiently strong means that the power $\alpha$ in the special flow rapresentation is close to $1$. One might hope that absolute continuity could hold for all powers $\alpha>1/2$, but this is out of reach with the current techniques.}.
%Singularity of the spectrum is in stark contrast with the recent result in \cite{FFK} on flows on tori with a degenerated singularity (or stopping point), which are shown to have absolutely continuous (and actually countable Lebesgue) spectrum. 

% We remark that stopping points or non-degenerate fixed points (including centers) are known to produce mixing \cite{Ko:mix} (at rates which are expected to be polynomial, see e.g.~\cite{Fa1}), while typical minimal locally Hamiltonian flows with  \emph{non-degenerate} saddles are  known \emph{not} to be mixing  by the work of Scheglov \cite{Sch:abs} for genus two and Ulcigrai \cite{Ul:abs} for any genus. 
At the heart of our proof of Theorem~\ref{singular}, on the other hand, is a strengthening of results on absence of mixing (in particular of the works \cite{FL03,FL05} and \cite{Sch:abs}). As already mentioned, we hope that the techniques of \cite{Ul:abs} might be pushed to allow to apply the singularity criterium for typical flows in the open set $\mathscr{U}_{min}$ of minimal, uniquely ergodic (weakly mixing) but not mixing locally Hamiltonian flows.

In the open set   $\mathscr{U}_{\neg min}$, which consists of flows with non-degenerate singularities that are not minimal, %but the locally Hamiltonian flows
 but have %but
 \emph{several} minimal components, the nature of the spectrum (for the restriction of a typical flow to a minimal component) is unclear. These flows are indeed mixing, but with sub-polynomial rate (see \cite{Rav:mix}, which provides logarithmic upper bounds) and it is not clear whether to expect singularity or  absolute continuity of the spectrum.

\subsection*{Acknowledgements}
The author is part of SwissMAP ({\it The Mathematics of Physics} National Centre for Compentence in Research) and is currently supported by a SNSF (Swiss National Science Foundation) Grant No.~200021\_ 188617/1. Both are acknowledged for their support.
\section{Bibliography}

%%%BIBLIO
\begingroup
\renewcommand{\section}[2]{}

\end{document}